\documentclass[12pt]{article}
\usepackage[dvips]{graphicx}
\usepackage{latexsym,amsmath,amsfonts,amscd}
\usepackage{multicol}
\usepackage{amsthm, amssymb,mathrsfs}
\usepackage{indentfirst} % indent the first sentence in each section
\usepackage{psfrag}

\newtheorem{theorem}{Theorem}[section]

\theoremstyle{remark}
  
\newtheorem{Example}{Example}[section]

\numberwithin{equation}{section}

\newcommand{\beq}{\begin{equation}}
\newcommand{\eeq}{\end{equation}}
\newcommand{\bed}{\begin{displaymath}}
\newcommand{\eed}{\end{displaymath}}
\newcommand{\bea}{\begin{array}}
\newcommand{\eea}{\end{array}}

\newcommand{\bld}[1]{\boldsymbol{#1}}
\newcommand{\veps}{\varepsilon}

\newcommand{\Oh}{{\Omega_h}}

\newcommand{\ol}[1]{\overline{#1}}                   % complex conjugate
\newcommand{\jmp}[1]{[\![#1]\!]}                     % jump
\newcommand{\bi}{\bld{i}}                            % imaginary unit
\usepackage{color}
\definecolor{red}{rgb}{1,0,0}

\definecolor{blck}{rgb}{0,0,0}

\theoremstyle{remark}

\begin{document}

\title{Numerical investigations on the resonance errors of  multiscale discontinuous Galerkin methods for one-dimensional stationary Schr\"{o}dinger equation
\\
%	Simulation of Quantum Transport
}

\author{Bo Dong\footnote{E-mail: bdong@umassd.edu. Department of Mathematics, University of Massachusetts Dartmouth, North Dartmouth, MA 02747.
%Research supported by NSF grant DMS-1818998.% Corresponding author: B. Dong.
}
and
Wei Wang\footnote{Corresponding author.
E-mail: weiwang1@fiu.edu. Department of Mathematics \& Statistics, Florida International University, Miami, FL, 33199.
} 
}

\date{}

\maketitle

\begin{abstract}
	%We carry out numerical experiments to investigate the impact of penalty parameters in the numerical traces on the resonance errors of multiscale discontinuous Galerkin (DG) methods  \cite{BW16,BW20} for one-dimensional Schr\"odinger equation. 
	%Using zero penalty parameters can simplify the implementation of the multiscale DG method, but it leads to obvious resonance errors when the mesh size $h$ is at a similar scale to the wave length. 
	%\bo{In previous work,  penalty parameters were assumed to be positive in error analysis, and numerical tests showed that the  methods also worked for zero penalty parameters. In this paper,  our numerical experiments  indicate that zero penalty parameter may lead to resonance errors for high order multiscale DG methods, but} taking positive penalty parameters can effectively reduce resonance errors and make the matrix in the global linear system have \ww{better} condition numbers.\\
%	
     In this paper, numerical experiments are carried out to investigate the impact of penalty parameters in the numerical traces on the resonance errors of high order multiscale discontinuous Galerkin (DG) methods  \cite{BW16,BW20} for one-dimensional stationary Schr\"odinger equation. 	
	 Previous work showed that penalty parameters were required to be positive in error analysis, but the  methods with zero penalty parameters worked fine in numerical simulations  on coarse meshes. In this work,  by performing extensive numerical experiments, we discover that  zero penalty parameters lead to resonance errors in the multiscale DG methods, and taking positive penalty parameters can effectively reduce resonance errors and make the matrix in the global linear system have {better} condition numbers. \\

\textbf{Key words:} discontinuous Galerkin method, multiscale method, resonance errors, one-dimensional Schr\"odinger equation
\end{abstract}

%%%%%%%%%%%%%%%%%%%%%%%%%%%%%%%%%%%%%%%%%%%%%%%%%%%%%%%%%%%%%%%%%%%%%%%%%%%%%%%%
\section{Introduction}
In this paper, we consider
the following one-dimensional second-order 
equation
\begin{equation}\label{eq:main}
-\varepsilon^2u''-f(x)u=0,
\end{equation}
where $\varepsilon>0$ is a small parameter and $f(x)$
is a real-valued smooth function.
One example of this type of equation is the
stationary Schr\"{o}dinger equation
in the modeling of quantum transport in nanoscale semiconductors  \cite{AbPi06,Ne08,WangShu09}.
%\begin{equation}\label{eq:schrod1d}
%\left\{
%\begin{array}{l}
%-\frac{\hbar^2}{2m}\varphi''(x)-qV(x)\varphi(x)=E\varphi(x)\quad\;\mbox{on} \; [a,b]\\
%\hbar \varphi'(a)+\bi p(a)\varphi_p(a)=2\bi p(a), \quad \hbar \varphi_p'(b)=\bi p(b)\varphi(b),
%\end{array}
%\right.
%\end{equation}
%where $\hbar$ is the reduced Plank constant,
%the constant $m$ is the effective mass,
%$q$ is the elementary positive charge of the electron,
%$V(x)$ is the total electrostatic potential in the device, $E$ is the injection energy,
%$p(x)=\sqrt{2m(E+qV(x))}$ is the momentum, and the solution $\varphi$ is a wave function.
%By letting $\displaystyle \varepsilon=\frac{\hbar}{\sqrt{2mE}}$ and $\displaystyle f(x)=1+\frac{qV(x)}{E}$, the Schr\"{o}dinger equation \eqref{eq:schrod1d} can be written into the model equation \eqref{eq:main}.

Note that when $f$ is positive, the solution to Eq. \eqref{eq:main} is an oscillatory wave function and the wave length is at the scale of $\varepsilon$. For small $\varepsilon$, the solution is highly oscillatory and simulation using standard finite element or finite difference methods requires very fine meshes. 
%In many applications, such \ww{an} equation needs to be solved for many times for different $\varepsilon$ and thus it is desirable to have numerical methods that can effectively approximate the wave solutions even on coarse meshes. 
{In applications such as quantum transport, this equation needs to be repeatly solved for  different magnitudes of wave lengths and thus it is desirable to develop multiscale methods that can effectively approximate the wave solutions even on coarse meshes. }

%---------
%Schr\"{o}dinger equation is frequently used in the modeling of quantum transport in nanoscale semiconductors  \cite{AbPi06,Ne08,WangShu09}.
%Schr\"{o}dinger equation involves small scales of wavelengths causes numerical challenges.  The oscillatory wave functions require very fine meshes to resolve for traditional numerical methods.

%----------

There have been ongoing effects on developing multiscale finite element methods which perform better than their polynomial-based counterparts for equations with highly oscillatory solutions; see \cite{Aarnes05,AbMN07,Gabard07,Perugia09,EfendievHou09, NBPM04,PoAb05,Yuan06, YS2, Wang11,Wang14} and references therein. 
 One of the popular ideas is to  incorporate the  information of small scale into  non-polynomial basis functions so that oscillatory solutions can be captured on coarse meshes. In \cite{AAN11,Ne08}, a second-order continuous finite element method based on WKB asymptotics was introduced and analyzed for solving the stationary Schr\"odinger equation. Similar multiscale methods in the discontinuous Galerkin (DG) framework were developed and analyzed for Eq. \eqref{eq:main} in \cite{WangShu09, BW16}. In \cite{BW20}, the multiscale DG methods were extended to higher orders using two different approximate function spaces $E^p$  and $T^{2p-1}$ for $p\ge 1$. Thanks to the advantage of having no continuity constraints across element interfaces, the multiscale DG methods in   \cite{WangShu09,BW20} have been extended to two dimensions in   \cite{GuoXu14,BW22,BW22proceeding}.  

In the multiscale DG methods, the numerical traces have been written in the form of alternating fluxes with penalty terms. 
%On one hand, the error analysis in  \cite{BW16, BW20, BW22, BW22proceeding} requires that the penalty parameters be nonzero (\ww{positive?}) for the proof of error estimates. On the other hand, numerical results show that the methods still work if the penalty parameters are zero.
 In \cite{WangShu09}, the third-order multiscale DG method using the $E^2$ space with zero penalty parameters %showed 
  produced good results in simulation of  Schr\"odinger equations in Resonant Tunneling Diode  on coarse meshes. However, 
the error analysis in later work \cite{BW16, BW20, BW22, BW22proceeding} %requires 
 required that the penalty parameters be positive  for proving the error estimates.  Furthermore, % (and the implementation of the method in such case will be simpler).
  numerical results in \cite{BW16} showed
 that the multiscale DG method using the space $E^2$ has resonance errors  around $h\simeq \varepsilon$ when the penalty parameters are zero. 
 
 In this paper, we perform numerical experiments to further    investigate  how the  values of penalty parameters affect the accuracy of general  multiscale DG methods using $E^p$ and $T^{2p-1}$ spaces.   Our numerical experiments indicate that if the exact solution lies in the approximation space, the approximate solutions from different choices of penalty parameter are all accurate up to rounding off errors.  But in the general case when the exact solution is not in the approximation space, zero penalty parameter will result in noticeable resonance errors when $h\simeq\varepsilon$, while using positive penalty parameters will greatly reduce such resonance errors. This improvement is also observed   in the condition number of the matrix in the resulting global linear system.

The rest of the paper is organized as follows. In section \ref{sec:method}, we describe the multiscale DG methods and related error estimates in previous work. In section \ref{sec:numerical_tests}, we demonstrate numerical results for  different choices of penalty parameters. The concluding remarks and future work are in Section 4.

%%%%%%%%%%%%%%%%%%%%%%%%%%%%%%%%%%%%%%%%%%%%%%%%%%%%%%%%%%%%%%%%%%%%%%%%%%%%%
\section{Multiscale DG method}\label{sec:method}

\subsection{The DG formulation}

We consider the equation \eqref{eq:main} with open boundary conditions, 
\begin{equation}\label{eq:u}
\left\{
\begin{array}{l}
-\varepsilon^2u''-f(x)\,u=0, \quad x\in [a,b],\\
%u'(a)+\bi k(a) \, u(a)= 2\bi k(a), \quad u'(b)-\bi k(b)\, u(b)=0,
\varepsilon u'(a)+\bi \sqrt{f_a}\;u(a)=2\bi \sqrt{f_a},\quad 
\varepsilon u'(b)-\bi\sqrt{f_b}\; u(b)=0,
\end{array}
\right.
\end{equation}
where $f_a=f(a)$, $f_b=f(b)$, and $\bi$ is the imaginary unit.
%where
%$\displaystyle k(x)=\sqrt{f(x)}/\varepsilon$ is the wave number.

By introducing an auxiliary variable $q$, we can rewrite the second-order equation in \eqref{eq:u}
 into the mixed form
\begin{subequations}\label{eq:mixed_form}
\begin{alignat}{1}
\label{eq:dgscheme}
&q-\varepsilon u'=0,\quad -\varepsilon q'-f(x)u=0,\\
\intertext{and the boundary conditions as}
\label{eq:bc2}
&q(a)+\bi \sqrt{f_a}\, u(a)=2\bi \sqrt{f_a},\quad q(b)-\bi \sqrt{f_b} \, u(b)=0.
\end{alignat}
\end{subequations}
%where $f_a=f(a)$ and $f_b=f(b)$.

To define multiscale DG methods for \eqref{eq:mixed_form},  we  need to introduce some notation. 
Let $a=x_{\frac{1}{2}}< x_{1+\frac{1}{2}}<\dots<x_{N+\frac{1}{2}}=b$ be a partition of the domain, $\Oh:=\{ I_j=(x_{j-\frac{1}{2}},x_{j + \frac{1}{2}}): j=1,\dots,N\}$ be  the set of all elements, 
%$\dOh:=\{ \partial I_j: j=1,\dots,N\}$   the set of boundaries of all elements, $\Eh:=\{x_{j+\frac{1}{2}}\}_{j=0}^N$ the set of all the \bo{element} interfaces, and $\Ehi:=\{x_{j+\frac{1}{2}}\}_{j=1}^{N-1}$ the set of all interior \bo{element} interfaces.  We use
%$x_j:=\frac{1}{2}(x_{j-\frac{1}{2}}+x_{j+\frac{1}{2}})$  the middle point of each element $I_j$, and
%$h =\max\limits_{j=1,\cdots, N} x_{j+\frac{1}{2}}-x_{j-\frac{1}{2}}$ the mesh size.
%$h_j=x_{j+\frac{1}{2}}-x_{j-\frac{1}{2}}$, 
%and $h =\max\limits_{j=1,\cdots, N} h_j$.
and  $h =\max\limits_{j=1,\cdots, N} (x_{j+\frac{1}{2}}-x_{j-\frac{1}{2}})$.

We consider two families of multiscale finite element spaces that we introduced in \cite{BW20}. 
The first approximate function space $E^p$ is obtained by combining two exponential basis functions and polynomial basis. For any element $I_j\in\Omega_h$, we let
\begin{subequations}
	\begin{alignat*}{1}
		E^{p}|_{I_j}=\begin{cases}
			\mbox{span}\{e^{\bi k_j (x-x_j)}, e^{-\bi k_j (x-x_j)}\}	& \mbox{ if } p=1,\\
			\mbox{span}\{  e^{\pm \bi k_j(x-x_j)}, 1, x, \cdots, x^{p-2} \} & \mbox{ if } p\ge 2,
		\end{cases}
	\end{alignat*}
\end{subequations}
where $\displaystyle k_j=\sqrt{f(x_j)}/\varepsilon$ is the wave number, and $x_j:=\frac{1}{2}(x_{j-\frac{1}{2}}+x_{j+\frac{1}{2}})$ is the middle point of the element $I_j$.
The second approximate function space $T^{2p-1}$ contains exponential basis functions in pairs. On each element $I_j\in\Omega_h$,
\begin{subequations}
	\begin{alignat*}{1}
		{T}^{2p-1}|_{I_j}=  \mbox{span}\{e^{\pm\bi k_j (x-x_j)},
		e^{\pm 2\bi k_j (x-x_j)},
		\cdots,
		e^{\pm p \bi k_j (x-x_j)}&\}
	\end{alignat*}
\end{subequations}
for any $p\ge 1$. 
Note that $T^1=E^1$, and the basis functions of $E^p$ and $T^{2p-1}$ are globally discontinuous across element interfaces and contains the information on the small scale of the problem.

Our DG method for \eqref{eq:mixed_form} is to seek approximate solutions $u_h, q_h\in E^p$ or $T^{2p-1}$   that satisfy the following weak formulation 
\begin{equation}
\label{eq:weak_form}
\begin{aligned} 
\sum_{j=1}^N\int_{I_j} q_h\,\ol{w}\,dx+\veps\sum_{j=1}^N\int_{I_j}\, u_h\,\ol{w'} \,dx-
\veps\sum_{j=1}^N\, \widehat{u}_h\,\ol{w}\,\bigg|^{x_{j+\frac{1}{2}}}_{x_{j-\frac{1}{2}}}&=0,\\
\veps\sum_{j=1}^N\int_{I_j} q_h\,\ol{v'}\,dx
-\veps\sum_{j=1}^N\, \widehat{q}_h\,\ol{v}\,\bigg|^{x_{j+\frac{1}{2}}}_{x_{j-\frac{1}{2}}}
-\sum_{j=1}^N\int_{I_j}\, f\,u_h\,\ol{v} \,dx&=0 
%(q_h, w)_\Oh+ (\veps u_h, w')_\Oh-\langle \veps \widehat{u}_h,w\,\bn\rangle_\dOh&=0\\
%(\veps q_h, v')_\Oh-\langle \veps \widehat{q}_h, v\,\bn\rangle_\dOh-(f(x)u_h, v)_\Oh&=0.
\end{aligned}
\end{equation}
for all test functions $v_h, w_h \in E^p$ or $T^{2p-1}$, 
%where we have used
%the notation
%\begin{alignat*}{1}
% (\varphi, v)_{I_j}=&\int_{I_j} \varphi(x)\ol{v}(x)\,dx, \\
%\langle \psi, w\,\bn \rangle_{\partial I_j}=&\psi(x_{j+\frac{1}{2}})\ol{w}(x_{j+\frac{1}{2}})-\psi(x_{j-\frac{1}{2}})\ol{w}(x_{j-\frac{1}{2}}),\\
%(\varphi, v)_{\Oh} =&\sum_{j=1}^N (\varphi, v)_{I_j},\quad
%\langle \psi, w\,\bn \rangle_{\dOh} = \sum_{j=1}^N \langle \psi, w\,\bn \rangle_{\partial I_j},
% \end{alignat*}
where  $\ol{w}$ is the complex conjugate of $w$.
% and $\bn$ is the unit outward normal vector. For  $I_j=(x_{j-\frac{1}{2}}, x_{j+\frac{1}{2}})$, we assume $\bn(x_{j-\frac{1}{2}})=-1$ and $\bn(x_{j+\frac{1}{2}})=1$.
%The numerical traces $\widehat{u}_h$ and $\widehat{q}_h$ will be defined in the next subsection.

%
%The finite element space $V_h$ contains functions which
%are discontinuous across cell interfaces. For standard DG methods, these functions are piece-
%wise polynomials. For our multiscale DG method, the basis functions are constructed to be non-polynomial functions which incorporate the
%small scales to better approximate the oscillating solutions.
%The proposed high order multiscale finite element spaces $V_h$ will be defined in Section \ref{sec:basis}.
%Since the functions in $V_h$ are allowed to have discontinuities at cell interfaces
%$x_{j+\frac{1}{2}}$,
%we use
% $v^-(x_{j+\frac{1}{2}})$ and  $v^+(x_{j+\frac{1}{2}})$ to refer to the left and right limits
%of $v$ at $x_{j+\frac{1}{2}}$, respectively.

%\subsection{The numerical traces}
%The choices of numerical traces are essential for the definition of DG methods, and different numerical traces will lead to different DG methods \cite{ABCM02}.
% In our schemes, we  use the same numerical traces as in our previous work \cite{BW18}. At the interior element interfaces,
The numerical traces $\widehat{u}_h$ and $\widehat{q}_h$ are defined in the same way as in \cite{BW16, BW20}, that is, at the interior element interfaces,
\begin{equation}\label{eq:traces_interior}
\begin{aligned} 
\widehat{u}_h(x_{j+\frac{1}{2}}) =\, &u_h^- (x_{j+\frac{1}{2}})- \bi \,\beta \,\jmp{q_h}(x_{j+\frac{1}{2}}),\\
\widehat{q}_h(x_{j+\frac{1}{2}})=\, &q_h^+(x_{j+\frac{1}{2}}) + \bi \,\alpha \,\jmp{u_h}(x_{j+\frac{1}{2}}),
\end{aligned}
\end{equation}
where
$v^-(x_{j+\frac{1}{2}})$ and  $v^+(x_{j+\frac{1}{2}})$ are the left and right limits
of $v$ at $x_{j+\frac{1}{2}}$, respectively, and $\jmp{v}=v^--v^+$ represents the jump across the interface. 
%We take the penalty parameters $\beta$ and $\alpha$ to be positive constants, which makes the DG methods different from the multiscale local DG in \cite{WangShu09}. This allows us to carry out error analysis in a way similar to  \cite{FengWu09}. Moreover, our numerical tests show that  the nonzero penalty terms are necessary for %non-resonance errors.
%\bo{reducing resonance errors.}

At the two boundary points $\{a, b\}$,
\begin{equation}\label{eq:traces_boundary}
\begin{aligned}
\widehat{u}_h(a)&=(1-\gamma)u_h(a)+ \bi \frac{\gamma}{\sqrt{f_a}}\, q_h(a) + 2\gamma,\\
\widehat{q}_h(a)&=\gamma q_h(a)-\bi (1-\gamma)\sqrt{f_a}\, u_h(a)+2\bi\,(1-\gamma)\sqrt{f_a},\\
\widehat{u}_h(b)&=(1-\gamma)u_h(b)- \bi \frac{\gamma}{\sqrt{f_b}}\, q_h(b),\\
\widehat{q}_h(b)
&=\gamma q_h(b)+\bi (1-\gamma)\sqrt{f_b}\, u_h(b),
\end{aligned}
\end{equation}
where $\gamma$ can be any real constant in $(0,1)$. %The numerical traces at the boundary are defined in this way to match the boundary conditions in \eqref{eq:bc2}.

Note that when the penalty parameters $\alpha$ and $\beta$ in \eqref{eq:traces_interior} are taken to be zero, the numerical traces at interior interfaces are alternating fluxes and
  the method reduces to MD-LDG method \cite{CoDo07} with multiscale bases.
  Note also that the numerical traces in \eqref{eq:traces_boundary} satisfy
$$ \widehat{q}_h(a)+\bi\mbox{}\sqrt{f_a}\mbox{ }\widehat{u}_h(a)=2\bi\sqrt{f_a},\quad
 \widehat{q}_h(b)-\bi\mbox{}\sqrt{f_b}\mbox{ }\widehat{u}_h(b)=0,$$
which coincide with the boundary conditions  \eqref{eq:bc2} for the exact solution.

\subsection{Error estimates in previous work}
We have proved some error estimates for the multiscale DG method defined by \eqref{eq:weak_form}-\eqref{eq:traces_boundary} in our  previous work.
For ease of reading, we list them here.
First, in \cite{BW16} we proved the following error estimates for the multiscale DG method with the approximation space $E^1$ .
\begin{theorem}\label{thm:E1}
	Let $u$ be the solution of the problem \eqref{eq:u} and $u_h\in E^1$ be the multiscale DG approximation defined by \eqref{eq:weak_form}-\eqref{eq:traces_boundary}.
	Assume that $\alpha$ and $\beta$ are positive constants and $0<\gamma<1$. For any mesh size $h>0$, we have
	$$\|u-u_h\| \le C |f|_{1,\infty}(\frac{h^2}{\veps}+\frac{h^3}{\veps^2})\|u\|,$$
	where $C$ is independent of $\veps$ and $h$.
\end{theorem}
This theorem  shows that when $f$ is constant, the approximate solution from the multiscale DG method will only have round-off errors. When $f$ is not constant, the method has a second order convergence rate for any mesh size, including $h\gtrsim \veps$.

For the multiscale DG method with higher order approximation spaces, we obtained the following error estimates in \cite{BW20}.
\begin{theorem}\label{thm:1Dhigh}
Suppose $(u, q)$ is the exact solution of the problem \eqref{eq:mixed_form}, $(u_h, q_h)$ is the solution of the multiscale DG method using $E^p$ for $p=2$ or 3, and $(\tilde{u}_h, \tilde{q}_h)$ is the solution of the multiscale DG method using $T^{2p-1}$ for $p=2$ or 3.  Assume that $\alpha$ and $\beta$ are positive constants and $0<\gamma<1$. When $h$ is small enough, we have
$$\|u-u_h\| \le C \,h^{\min\{s+1, p+1\}}((\|u\|_{s+1}+\|q\|_{s+1}) ,$$
$$\|u-\tilde{u}_h\| \le C \,h^{\min\{s+1, 2p\}}((\|u\|_{s+1}+\|q\|_{s+1}),$$
	where $C$ is independent of $h$  but may depends on $\veps$.
\end{theorem}
Note that error estimates in Theorem \ref{thm:1Dhigh} are valid only if $h$ is sufficiently small. There are no error estimates for higher order spaces on coarse meshes.

The proofs of Theorem \ref{thm:E1} and Theorem \ref{thm:1Dhigh} require the penalty parameters   $\alpha$ and $\beta$ to be positive. % for technical reason. 
%
%Numerical tests in \cite{WangShu09} showed that the multiscale DG method using $E^2$  works for some numerical simulations if $\alpha$ and $\beta$ are taken as zero. In this case, the global matrix is banded and the method seems to be more efficient than using positive penalty parameters. 
%
Numerical simulations in \cite{WangShu09} showed that the multiscale DG method using $E^2$  with zero penalty parameters produced good results on coarse meshes.
In this case, the global matrix in the linear system is banded and the variable $q$ in Eq. \eqref{eq:dgscheme} can be solve locally, thus solving the linear system is more efficient than using positive penalty parameters. 
But in \cite{BW16}, we found that the multiscale DG method using $E^2$  with zero penalty parameters has obvious resonance errors around $h\simeq \varepsilon$.
%  and such errors can be reduced by using positive penalty parameters. 
%\ww{whereas multiscale DG methods with positive penalty parameters do not.}
Therefore, in the next section, we would like to numerically investigate how penalty parameters affect the errors of  multiscale DG methods with $E^p$ and $T^{2p-1}$ spaces, especially in the region $h\simeq \veps$.

\section{Numerical results}\label{sec:numerical_tests}

%\bo{Remove(}Recall in \cite{BW16}, we found that $E^2$ basis with zero jump coefficients has obvious resonance errors around $h\simeq \varepsilon$  and such errors can be reduced by using nonzero jump coefficients.\bo{)}
In this section, we carry out numerical experiments for 
multiscale DG 
%with different penalty parameters 
using   high order spaces  $E^{p}$ and $T^{2p-1}$, in particular,
 ${E}^1(=T^1)$,  ${E}^2$,  ${E}^3$,  ${T}^3$ and ${T}^5$.
  We will test the methods with different penalty parameters, compare the behaviors of the approximate solutions, and inspect the occurrence of resonance errors. 
We will also investigate the relation between the  resonance error and the condition number of the matrix resulting from the global linear system.

\Example \label{ex:1}

In the first example, we consider the simple case %with 
of constant
function $f(x)$.  
It is easy to see that in this case, the exact solution
of \eqref{eq:u} is in the   high-order finite element spaces  $E^{p}$ and $T^{2p-1}$. Thus the
multiscale DG method with these spaces can compute the solution exactly with
only round-off errors. 
We choose $f(x)=10$ which is the same as in our previous work \cite{BW16, BW20}. 
In \cite{BW20}, 
we performed the experiments of the multiscale DG methods only with 
 non-zero penalty parameters $\alpha=\beta=1$.
This time, we will compare the multiscale DG with three different magnitudes of  penalty parameters $\alpha=\beta=0$,
$\alpha=\beta=0.1$ and $\alpha=\beta=1$. Note that $\alpha=\beta=0$ reduces the method to the standard MD-LDG  method with multiscale bases.

The $L^2$-errors of the multiscale DG  method with $E^1$, ${E}^2$, ${E}^3$, $T^{3}$ and $T^{5}$ are shown
in Tables \ref{table:E1}-\ref{table:T5} respectively.  In each table, we compute the %method 
 approximate solutions with three different magnitudes of  penalty parameters $\alpha=\beta=0$,
$\alpha=\beta=0.1$ and $\alpha=\beta=1$ for two different levels of $\varepsilon$%:
,  
$\varepsilon=0.005$ and $\varepsilon=0.001$.
These two values of $\varepsilon$ are chosen 
because during the mesh refinement, $h$ is changing from  $h> \varepsilon$,  $h\simeq \varepsilon$,  to  $h< \varepsilon$.
 When $h\simeq \varepsilon$, the error may suddenly increase due to resonance. %it is called the resonance.
 
 For all the results in Tables  \ref{table:E1}-\ref{table:T5},
 it is clear to see the round-off errors
in double precision.
For the same value of $\varepsilon$, the results are very similar for different penalty parameters. Because they all capture the solution exactly, there is no
resonance errors in this example.
We notice  that when $\varepsilon$ is smaller, the round-off error increases slightly. %a little bit. 
This is because we use more points
 in each element for the numerical integration of the exponential functions which accumulate round-off errors for smaller $\varepsilon$.

 \begin{table}
 	\small
     \centering
   \caption{Example \ref{ex:1}: $L^2$-errors  by multiscale DG with $E^1$ for $f(x)=10$. }
    \smallskip
  \begin{tabular}{c|ccc|ccc}    \hline\hline
 &      \multicolumn{3}{c|}{$\varepsilon=0.005$}& \multicolumn{3}{c}{$\varepsilon=0.001$ }\\\hline
$N$&$\alpha=\beta=0$ & $\alpha=\beta=0.5$ &$\alpha=\beta=1$&$\alpha=\beta=0$&$\alpha=\beta=0.5$&$\alpha=\beta=1$\\\hline
10& 4.00E-12   & 4.17E-12     &  4.62E-12      & 7.08E-11          &  7.63E-11          &    7.66E-11        \\
20& 3.95E-12   &  4.01E-12    &     4.36E-12   &    6.92E-11       &   7.19E-11         &    7.27E-11         \\
%30& 3.76E-12   & 3.82E-12     &  3.90E-12      &  6.97E-11         & 7.31E-11           &    7.03E-11         \\
40&  4.53E-12  &  4.56E-12    &   4.78E-12     &  7.19E-11          &   7.26E-11         &     7.42E-11        \\
80& 4.58E-12   &  4.63E-12     &  4.61E-12      &   7.21E-11        &  7.34E-11          &     7.31E-11        \\
160& 2.62E-12   &  2.65E-12     &  2.65E-12      &   6.12E-11         &    6.15E-11         &   6.23E-11          \\
200& 6.71E-12   &  6.75E-12    &  6.77E-12      &   8.24E-11        &     8.23E-11        &    8.31E-11         \\\hline

    \end{tabular}\label{table:E1}
    \end{table}

%E1 30; E2 40; E3 25; T3 40
%Ne=200, 500
 \begin{table}
     \centering
   \caption{Example \ref{ex:1}: $L^2$-errors  by multiscale DG with $E^2$ for $f(x)=10$. }
    \smallskip
  \begin{tabular}{c|ccc|ccc}    \hline\hline
 &      \multicolumn{3}{c|}{$\varepsilon=0.005$}& \multicolumn{3}{c}{$\varepsilon=0.001$ }\\\hline
$N$&$\alpha=\beta=0$ & $\alpha=\beta=0.5$ &$\alpha=\beta=1$&$\alpha=\beta=0$&$\alpha=\beta=0.5$&$\alpha=\beta=1$\\\hline
10& 4.00E-12   & 4.18E-12     &  4.66E-12      &  7.08E-11         &   7.67E-11         &  7.67E-11          \\
20& 3.95E-12   &  4.01E-12    &     4.39E-12   &  6.93E-11         &  7.11E-11          &   7.31E-11          \\
40&  4.54E-12  &  4.56E-12    &   4.68E-12     &   7.19E-11        & 7.27E-11           &   7.45E-11           \\
80& 4.57E-12   &  4.61E-12     &  4.62E-12      &   7.21E-11        & 7.30E-11            &  7.42E-11           \\
%140& 4.30E-12   &  4.31E-12     &  4.31E-12      &           &            &             \\
160& 2.60E-12   &  2.61E-12     &  2.61E-12      &   6.14E-11        &  6.12E-11           &  6.19E-11           \\
200& 6.75E-12   &  6.75E-12    &  6.77E-12      &  8.22E-11         & 8.32E-11           &  8.22E-11            \\\hline

    \end{tabular}%\label{table:E2}
    \end{table}

\begin{table}
     \centering
   \caption{Example \ref{ex:1}: $L^2$-errors  by multiscale DG with $E^3$ for $f(x)=10$. }
    \smallskip
  \begin{tabular}{c|ccc|ccc}    \hline\hline
 &      \multicolumn{3}{c|}{$\varepsilon=0.005$}& \multicolumn{3}{c}{$\varepsilon=0.001$ }\\\hline
$N$&$\alpha=\beta=0$ & $\alpha=\beta=0.5$ &$\alpha=\beta=1$&$\alpha=\beta=0$&$\alpha=\beta=0.5$&$\alpha=\beta=1$\\\hline
10& 6.49E-12   & 4.19E-12     &  4.66E-12      &  4.66E-09         &   7.74E-11         &  8.24E-09          \\
20& 3.95E-12   &  4.01E-12    &     4.15E-12   &  6.93E-11         &  1.52E-09          &   7.24E-11          \\
40&  4.54E-12  &  4.58E-12    &   4.68E-12     &   7.19E-11        & 7.28E-11           &   7.44E-11           \\
80& 4.56E-12   &  4.58E-12     &  4.60E-12      &   7.20E-11        & 7.28E-11            &  7.29E-11           \\
%140& 4.30E-12   &  4.31E-12     &  4.31E-12      &           &            &             \\
160& 2.61E-12   &  2.62E-12     &  2.61E-12      &   6.10E-11        &  6.18E-11           &  6.09E-11           \\
200& 6.71E-12   &  6.72E-12    &  6.77E-12      &  8.22E-11         & 8.27E-11           &  8.27E-11            \\\hline

    \end{tabular}%\label{table:E3}
    \end{table}

\begin{table}
     \centering
   \caption{Example \ref{ex:1}: $L^2$-errors  by multiscale DG with $T^3$ for $f(x)=10$. }
    \smallskip
  \begin{tabular}{c|ccc|ccc}    \hline\hline
 &      \multicolumn{3}{c|}{$\varepsilon=0.005$}& \multicolumn{3}{c}{$\varepsilon=0.001$ }\\\hline
$N$&$\alpha=\beta=0$ & $\alpha=\beta=0.5$ &$\alpha=\beta=1$&$\alpha=\beta=0$&$\alpha=\beta=0.5$&$\alpha=\beta=1$\\\hline
10& 4.01E-12   & 4.17E-12     &  4.66E-12      &  7.09E-11         &   7.62E-11         &  7.56E-11          \\
20& 3.95E-12   &  4.01E-12    &     4.37E-12   &  6.98E-11         &  7.19E-11        &   7.24E-11          \\
40&  4.55E-12  &  4.56E-12    &   4.76E-12     &   7.19E-11        & 7.59E-11           &   7.44E-11           \\
80& 4.83E-12   &  4.68E-12     &  4.62E-12      &   7.24E-11        & 7.36E-11            &  7.75E-11           \\
%140& 4.30E-12   &  4.31E-12     &  4.31E-12      &           &            &             \\
160& 2.62E-12   &  2.64E-12     &  2.65E-12      &   6.13E-11        &  6.16E-11           &  6.74E-11           \\
200& 6.72E-12   &  6.72E-12    &  6.37E-12      &  8.22E-11         & 8.23E-11           &  8.35E-11            \\\hline

    \end{tabular}%\label{table:T3}
    \end{table}

\begin{table}
     \centering
   \caption{Example \ref{ex:1}: $L^2$-errors  by multiscale DG with $T^5$ for $f(x)=10$. }
    \smallskip
  \begin{tabular}{c|ccc|ccc}    \hline\hline
 &      \multicolumn{3}{c|}{$\varepsilon=0.005$}& \multicolumn{3}{c}{$\varepsilon=0.001$ }\\\hline
$N$&$\alpha=\beta=0$ & $\alpha=\beta=0.5$ &$\alpha=\beta=1$&$\alpha=\beta=0$&$\alpha=\beta=0.5$&$\alpha=\beta=1$\\\hline
10& 4.01E-12   & 4.19E-12     &  4.62E-12      &  7.57E-11       &   8.66E-11         &  8.58E-11          \\
20& 3.93E-12   &  4.01E-12    &     4.39E-12   &  6.98E-11         &  7.20E-11          &   7.05E-11          \\
40&  4.55E-12  &  4.58E-12    &   4.75E-12     &   7.18E-11        & 7.39E-11           &   7.39E-11           \\
80& 4.51E-12   &  4.58E-12     &  4.65E-12      &   7.23E-11        &7.37E-11            &  7.53E-11           \\
%140& 4.30E-12   &  4.31E-12     &  4.31E-12      &           &            &             \\
160& 2.62E-12   &  2.62E-12     &  2.65E-12      &   6.14E-11        &  6.16E-11           &  6.82E-11           \\
200& 6.72E-12   &  6.72E-12    &  6.77E-12      &  8.31E-11         & 8.30E-11           &  8.33E-11            \\\hline

    \end{tabular}\label{table:T5}
    \end{table}

\Example \label{ex:2}

Next, we consider the model equation  \eqref{eq:u} with a smooth  
positive function $f(x)=\sin{x}+2$ on $[0,1]$.
There is no simple explicit formula for the exact solution in this case and the solution is  not in the finite element spaces.
The reference solutions are computed by the MD-LDG method \cite{CoDo07} using piecewise cubic polynomials on $N=500,000$ elements.
We again test three different magnitudes of  penalty parameters   $\alpha=\beta=0$,
$\alpha=\beta=0.1$ and $\alpha=\beta=1$ 
for five different order spaces   ${E}^1 (=T^1)$,  ${E}^2$,  ${E}^3$,  ${T}^3$ and ${T}^5$.
We consider two values of $\varepsilon$: $\varepsilon=0.005$ and 0.001.
The mesh is refined from $N=5$ to $N=640$. %For $\varepsilon=0.005$, the resonance occurs around $N=20$. And for $\varepsilon=0.005$, the resonance occurs around N=200.

\begin{figure}
  \centering
  \includegraphics[width=2.6in,height=2.3in]{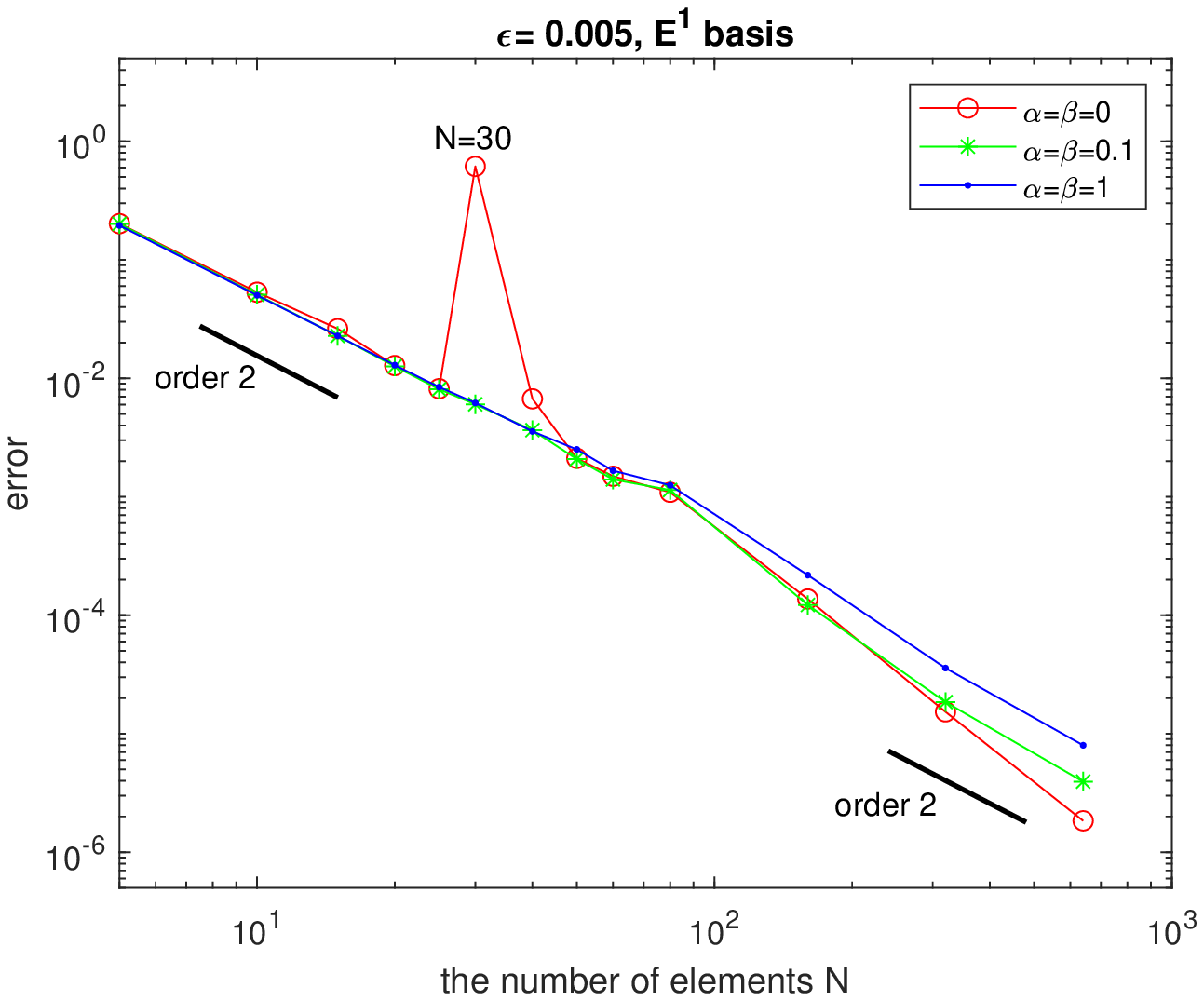}
    \includegraphics[width=2.6in,height=2.3in]{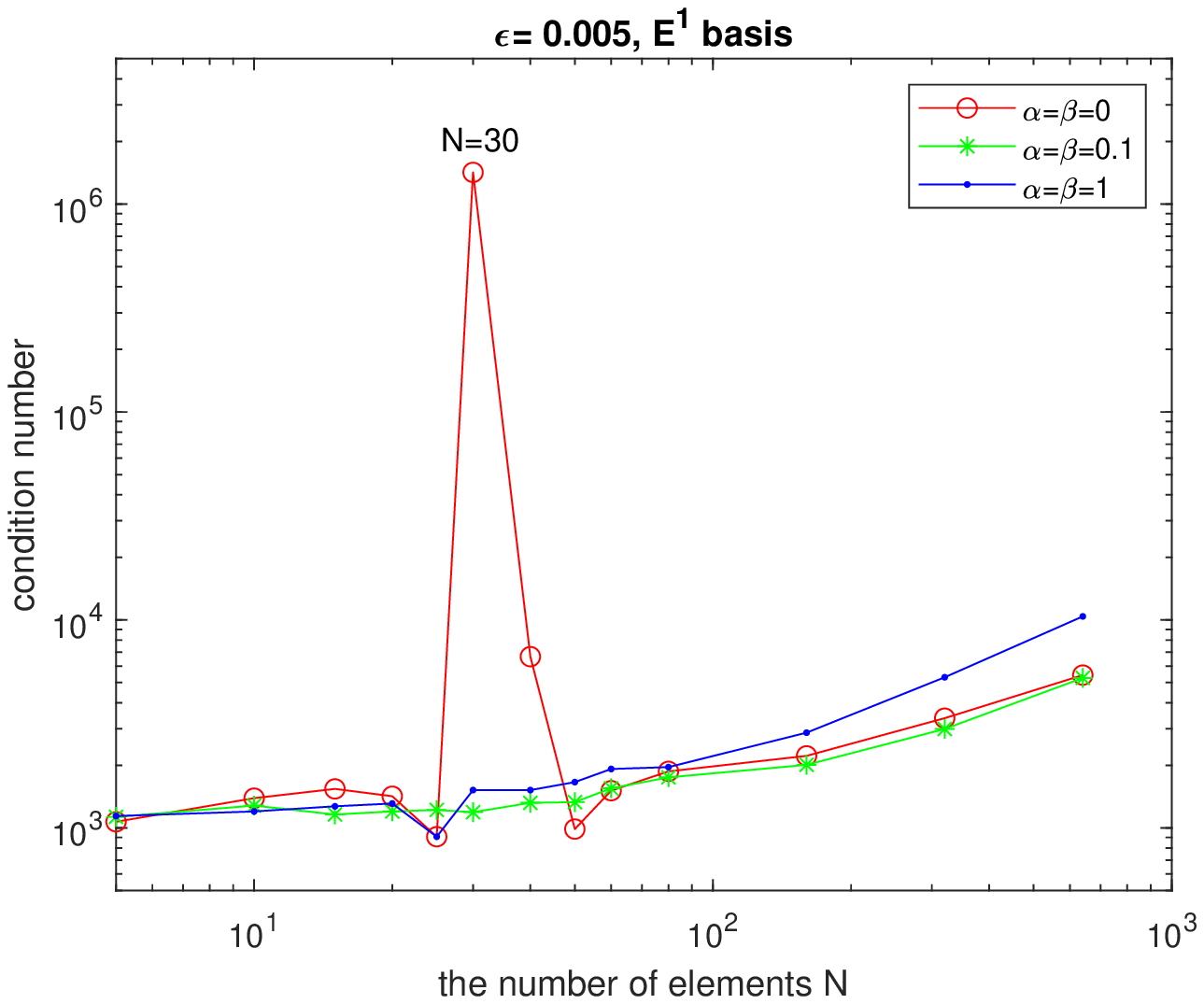}

\caption{Example \ref{ex:2}: Numerical results by multiscale DG ${E}^1$ for $\varepsilon=5\times 10^{-3}$. Left: $L^2$-errors of $u$. Right: condition numbers.}
\label{fig:5d3e1}
\end{figure}

Fig. \ref{fig:5d3e1} shows the log-log plot of  $L^2$-errors 
 %versus the number of elements $N$ for the multiscale DG using ${E}^1$  on the left and the corresponding condition numbers versus $N$ on the right.
   of $u$  versus the number of elements $N$  on the left and the corresponding condition numbers versus $N$ on the right for the multiscale DG  ${E}^1$ with $\veps=0.005$. In each figure, we compare the results from three different magnitudes of  penalty parameters $\alpha=\beta=0$, $\alpha=\beta=0.1$ and $\alpha=\beta=1$.
We can see the multiscale DG ${E}^1$ with  $\alpha=\beta=0$ (in red) has resonance errors around   $N=30$. Away from the resonance
location, 
%the multiscale DG ${E}^1$ 
the method shows a second order of convergence.
The multiscale DG ${E}^1$ with positive parameters  $\alpha=\beta=0.1$ (green) and $\alpha=\beta=1$ (blue)  do not show any resonance errors
  and both of them have %show 
the second order of convergence for all $N$, which is consistent with Theorem \ref{thm:E1}. 
Thus, we can see that positive penalty parameters significantly help reduce the resonance errors.

%We observe that the methods using $\alpha=\beta=0.1$ behave slightly better than $\alpha=\beta=1$ in terms of smaller errors and condition numbers.

Next, we look at the plot of condition numbers  on the right of Figure \ref{fig:5d3e1}. We see a big spike appearing at the same  
resonance location $N=30$ for the multiscale DG ${E}^1$ with  $\alpha=\beta=0$.  But the condition number plots for $\alpha=\beta=0.1$ and $\alpha=\beta=1$ are slowly increasing without spikes.
This suggests that resonance errors are related to the suddenly increased condition numbers.
We also observe that the method using $\alpha=\beta=0.1$ behaves slightly better than $\alpha=\beta=1$  for having smaller errors and condition numbers.

\begin{figure}
  \centering
  \includegraphics[width=2.6in,height=2.3in]{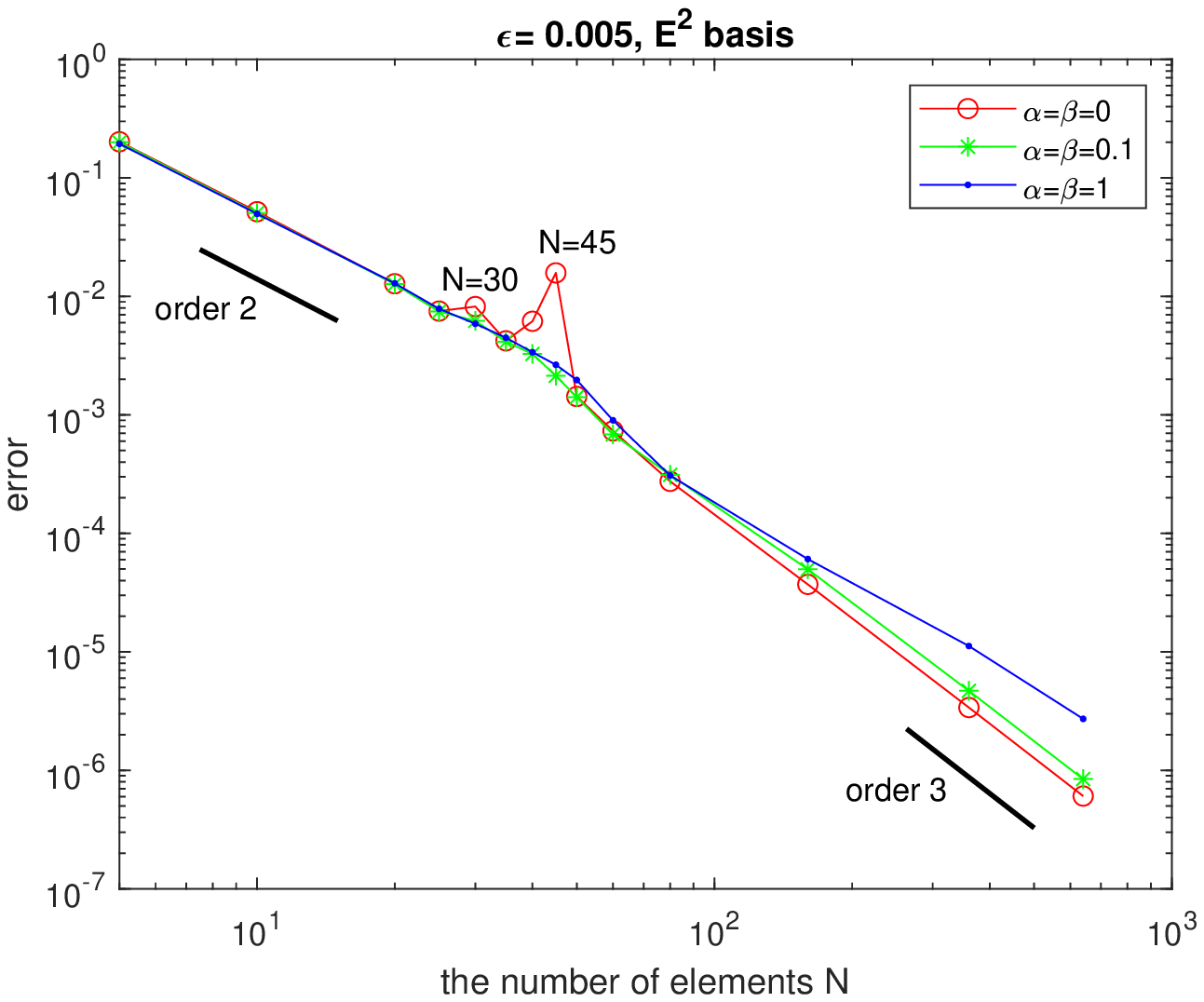}
    \includegraphics[width=2.6in,height=2.3in]{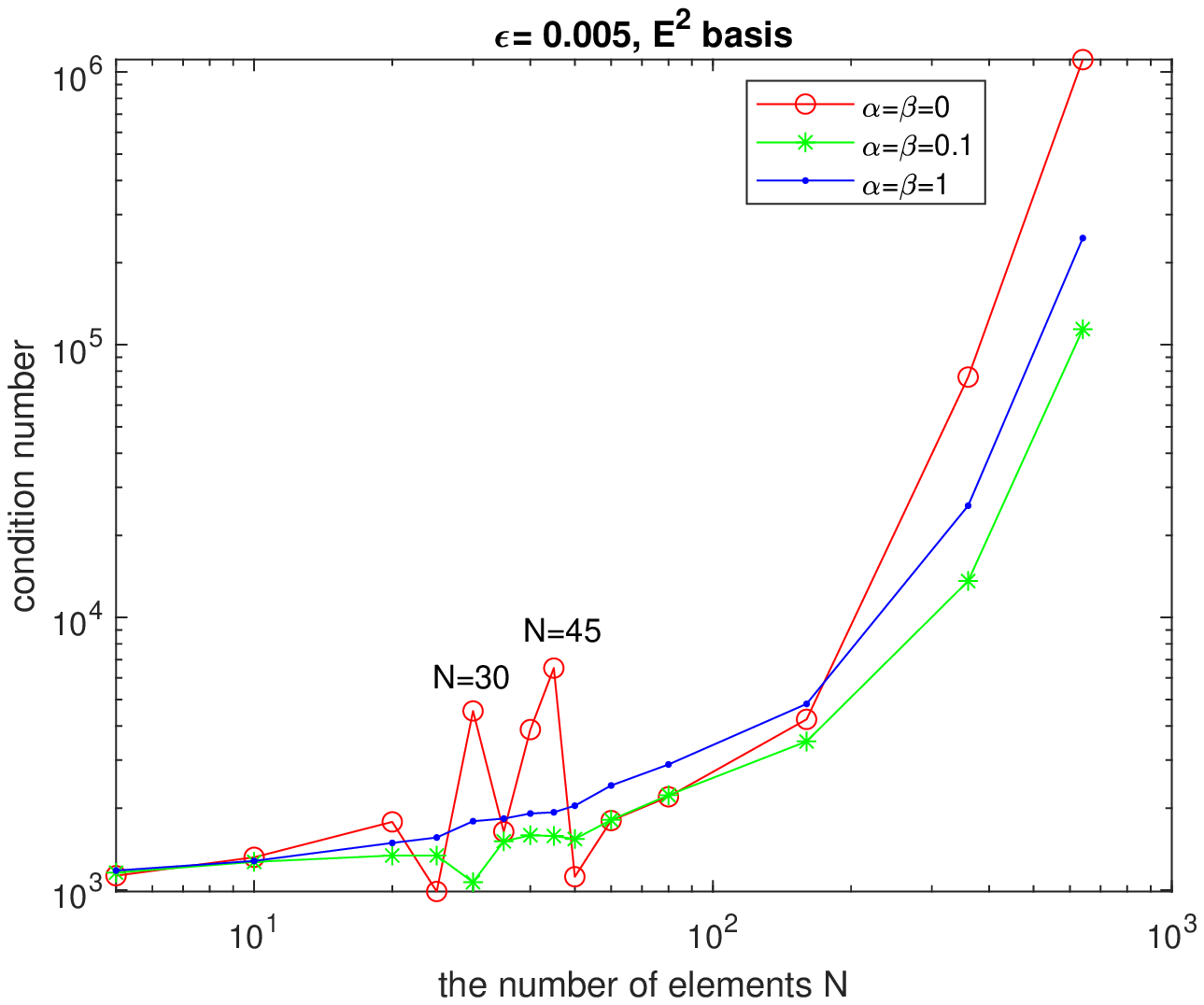}

\caption{Example \ref{ex:2}: Numerical results by multiscale DG ${E}^2$ for $\varepsilon=5\times 10^{-3}$. Left: $L^2$-errors of $u$. Right: condition numbers.}
\label{fig:5d3e2}
\end{figure}

We plot  %the same 
figures for multiscale DG ${E}^2$ in Fig.  \ref{fig:5d3e2} with $\varepsilon=0.005$.
We can see the method with $\alpha=\beta=0$ has resonance errors around   $N=30$ and  $N=45$ on the left of Fig.  \ref{fig:5d3e2}. At the same locations, 
condition number plot has spikes too (see Fig.  \ref{fig:5d3e2} right).
The multiscale DG ${E}^2$ with $\alpha=\beta=0.1$ and $\alpha=\beta=1$ do not have  resonance errors.
We observe a second order convergence when  $h\gtrsim \varepsilon$.
The approximate solutions show an optimal third order convergence when $h\lesssim \varepsilon$, which is consistent with Theorem \ref{thm:1Dhigh}. For zero penalty parameters, the method has similar convergence orders in the regions away from the resonance location.
%%multiscale DG ${E}^3$ in Fig.  \ref{fig:5d3e3},
%%multiscale DG ${T}^3$ in Fig.  \ref{fig:5d3t3}, multiscale DG ${T}^5$ in Fig.  \ref{fig:5d3t5}.

\begin{figure}
  \centering
  \includegraphics[width=2.6in,height=2.3in]{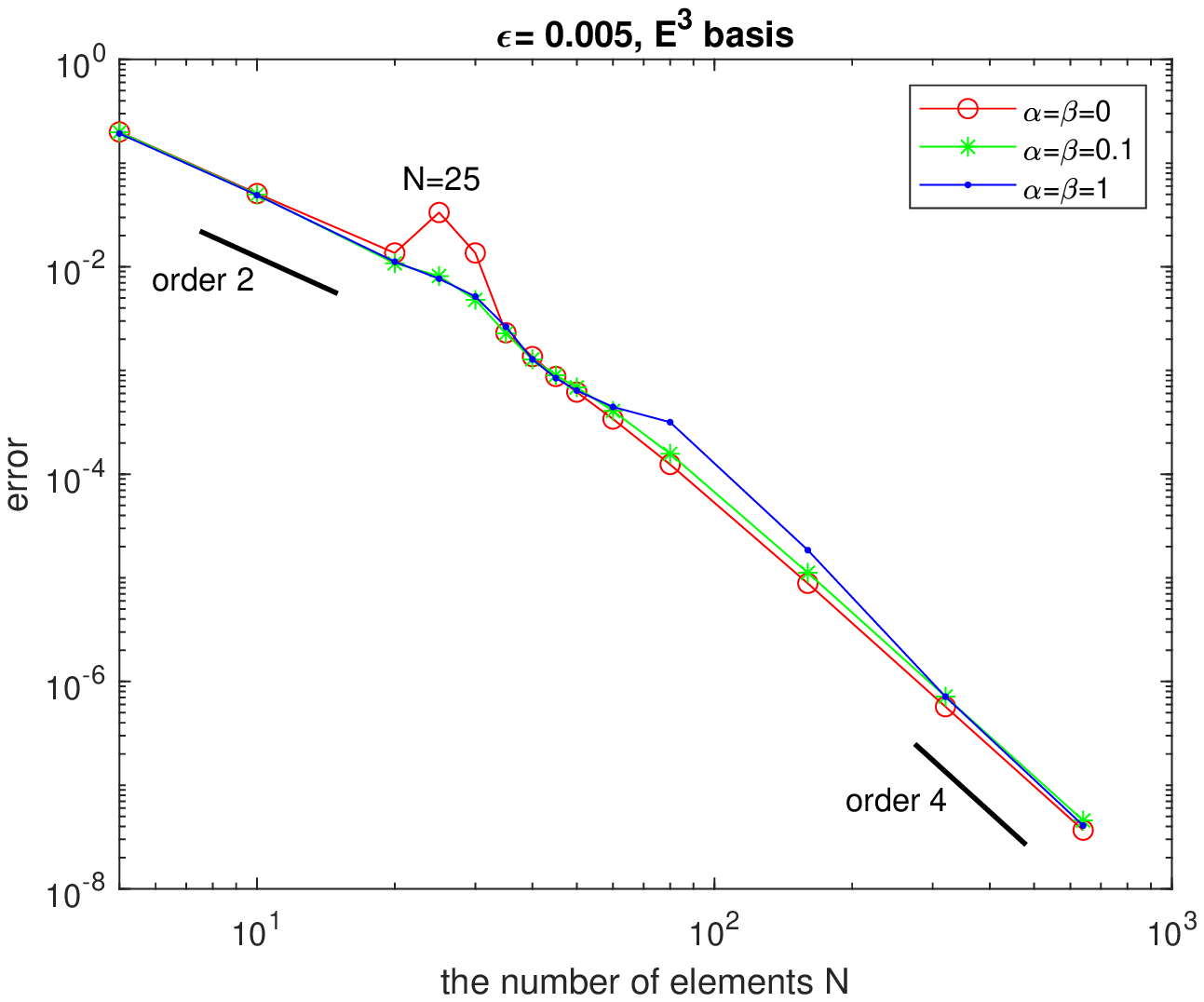}
    \includegraphics[width=2.6in,height=2.3in]{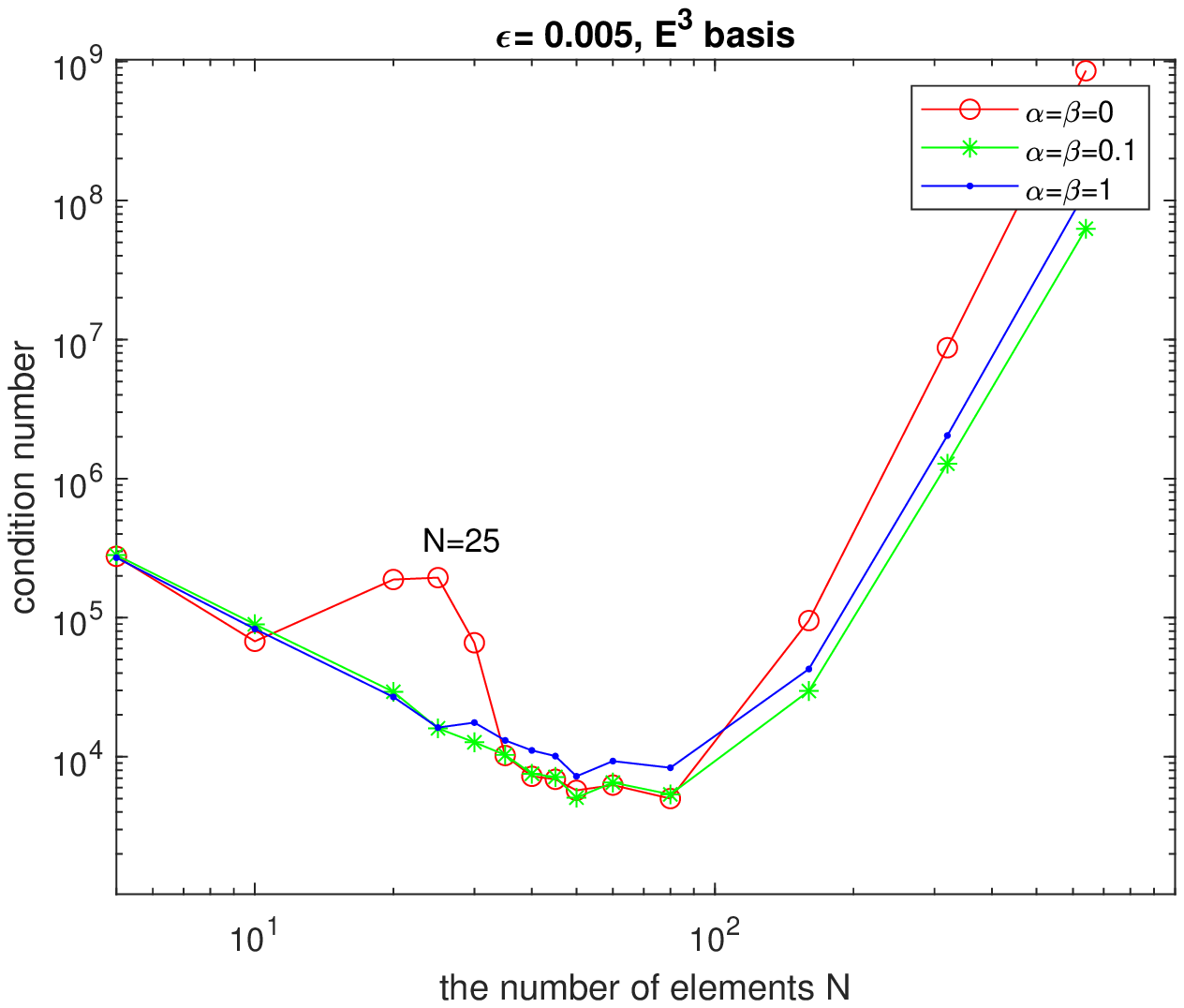}

\caption{Example \ref{ex:2}: Numerical results by multiscale DG ${E}^3$ for $\varepsilon=5\times 10^{-3}$. Left: $L^2$-errors of $u$. Right: condition numbers.}
\label{fig:5d3e3}
\end{figure}

We can see similar phenomenon for multiscale DG ${E}^3$ in Fig.  \ref{fig:5d3e3} with $\varepsilon=0.005$.
The method with $\alpha=\beta=0$ has resonance errors around   $N=25$ and  condition number has a big spike at the same location too.
Away from the resonance area, % $[25,30]$, 
we observe %(all multiscale DG ${E}^3$  have) 
 that all three choices of penalty parameters lead to a second order convergence when  $h \gtrsim \varepsilon$ and
 an optimal fourth order convergence when  $h \lesssim \varepsilon$.
%(We noticed a first increasing and then decreasing of condition numbers as $N$ increases in this case.) 
% We notice that the condition numbers first increase and then decrease as $N$ becomes larger.  We do not know the reason.

\begin{figure}
  \centering
  \includegraphics[width=2.6in,height=2.3in]{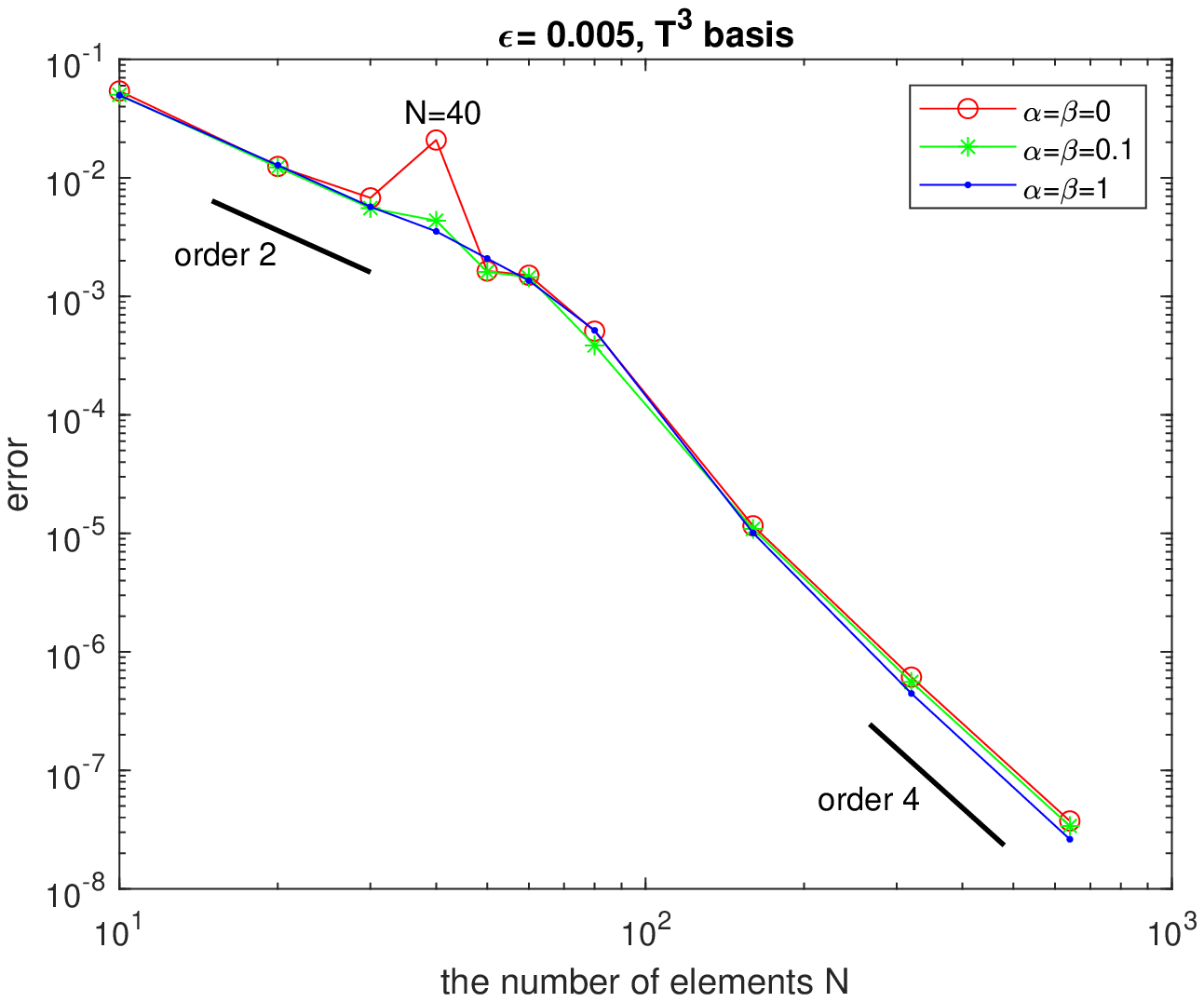}
    \includegraphics[width=2.6in,height=2.3in]{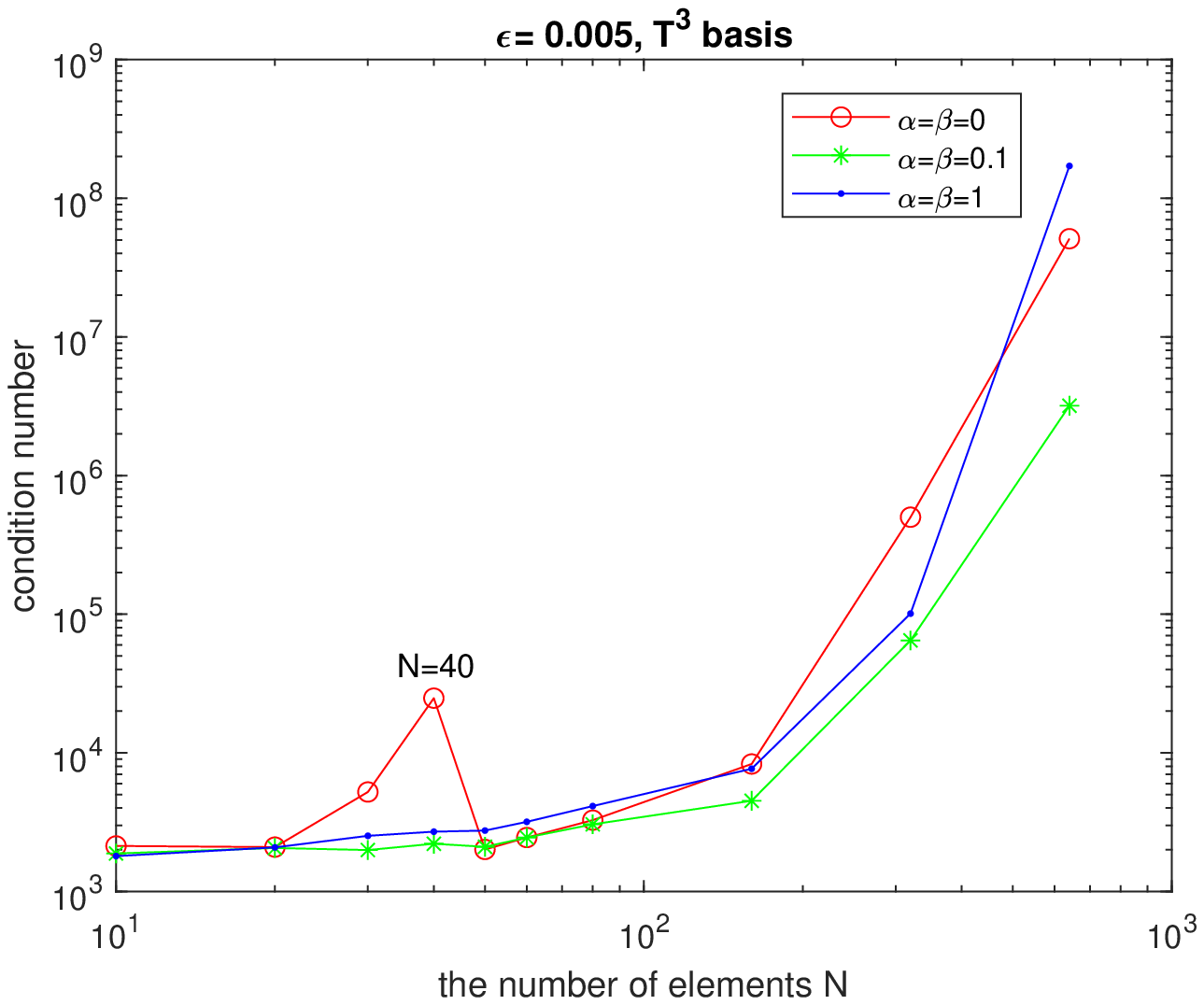}

\caption{Example \ref{ex:2}: Numerical results by multiscale DG ${T}^3$ for $\varepsilon=5\times 10^{-3}$. Left: $L^2$-errors of $u$. Right: condition numbers.}
\label{fig:5d3t3}
\end{figure}

For %(the same) 
 the other fourth order space ${T}^3$, the resonance errors occur around   $N=40$ in Fig.  \ref{fig:5d3t3}.
It shows a similar second order convergence when  $h \gtrsim \varepsilon$ and
 optimal fourth order convergence when  $h \lesssim  \varepsilon$.

\begin{figure}
  \centering
  \includegraphics[width=2.6in,height=2.3in]{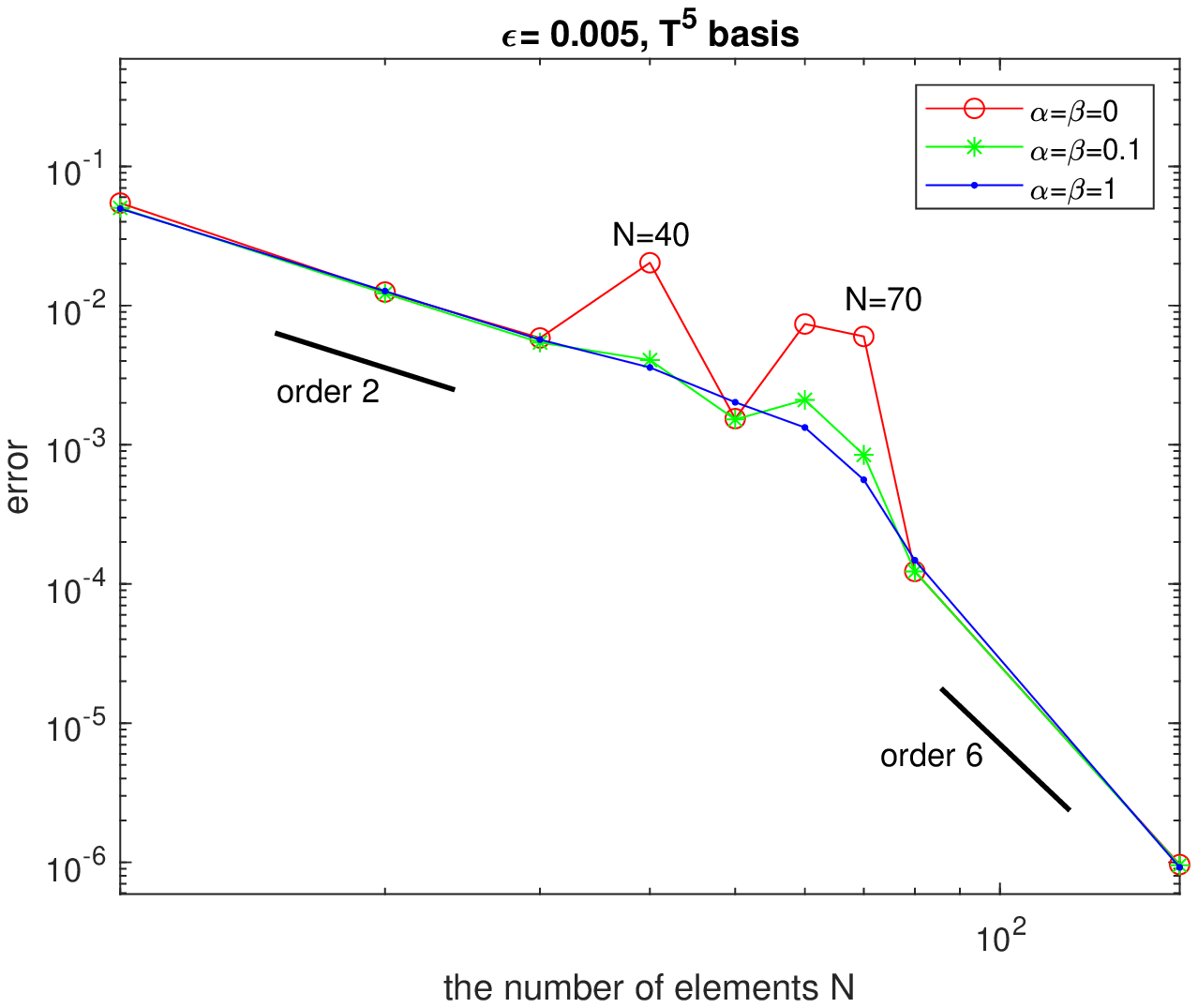}
    \includegraphics[width=2.6in,height=2.3in]{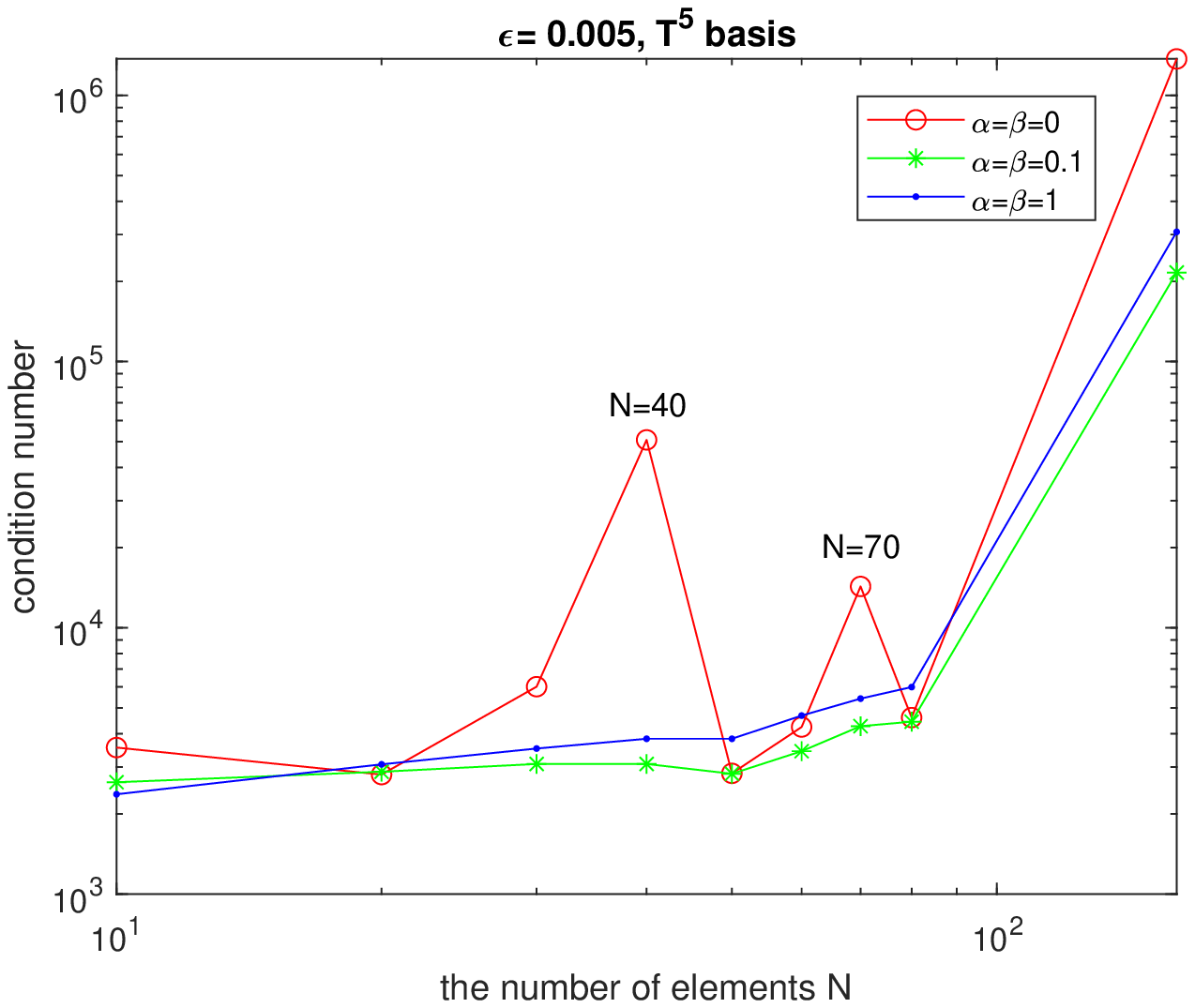}

\caption{Example \ref{ex:2}: Numerical results by multiscale DG ${T}^5$ for $\varepsilon=5\times 10^{-3}$. Left: $L^2$-errors of $u$. Right: condition numbers.}
\label{fig:5d3t5}
\end{figure}

 In Fig. \ref{fig:5d3t5}, the multiscale DG ${T}^5$ with $\alpha=\beta=0$ has two resonance errors regions around   $N=40$ and  $N=70$.
The condition number has   spikes at the same locations.
Away from the resonance area, 
we observe %(all multiscale DG ${T}^5$ show) 
 that the multiscale DG ${T}^5$ with all three choice of penalty parameters has a second order convergence when  $h \gtrsim \varepsilon$ and
  an optimal sixth order convergence when  $h \lesssim \varepsilon$.

\begin{figure}
  \centering
  \includegraphics[width=2.6in,height=2.3in]{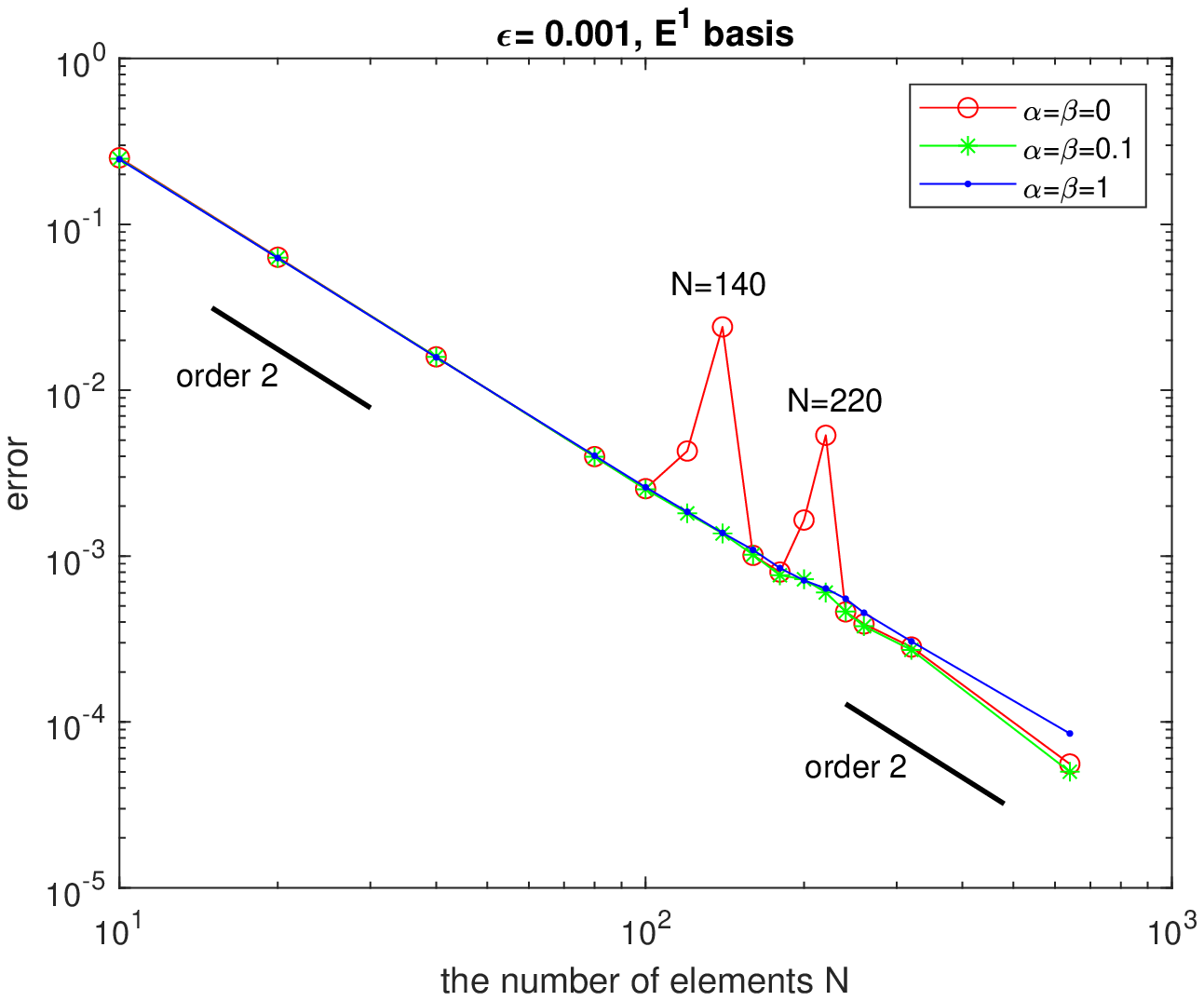}
    \includegraphics[width=2.6in,height=2.3in]{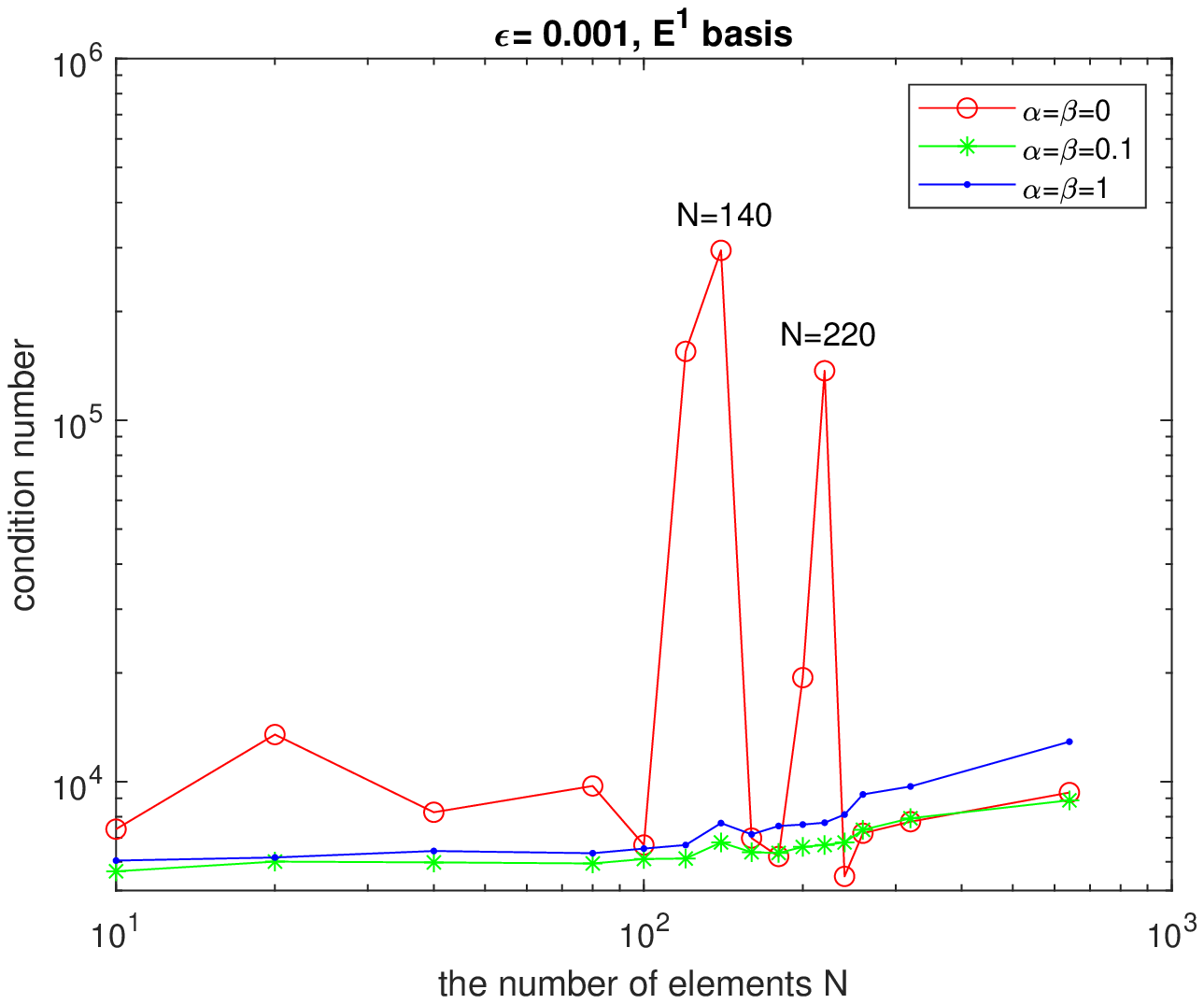}

\caption{Example \ref{ex:2}: Numerical results by multiscale DG ${E}^1$ for $\varepsilon=1\times 10^{-3}$. Left: $L^2$-errors of $u$. Right: condition numbers.}
\label{fig:1d3e1}
\end{figure}

\begin{figure}
  \centering
  \includegraphics[width=2.6in,height=2.3in]{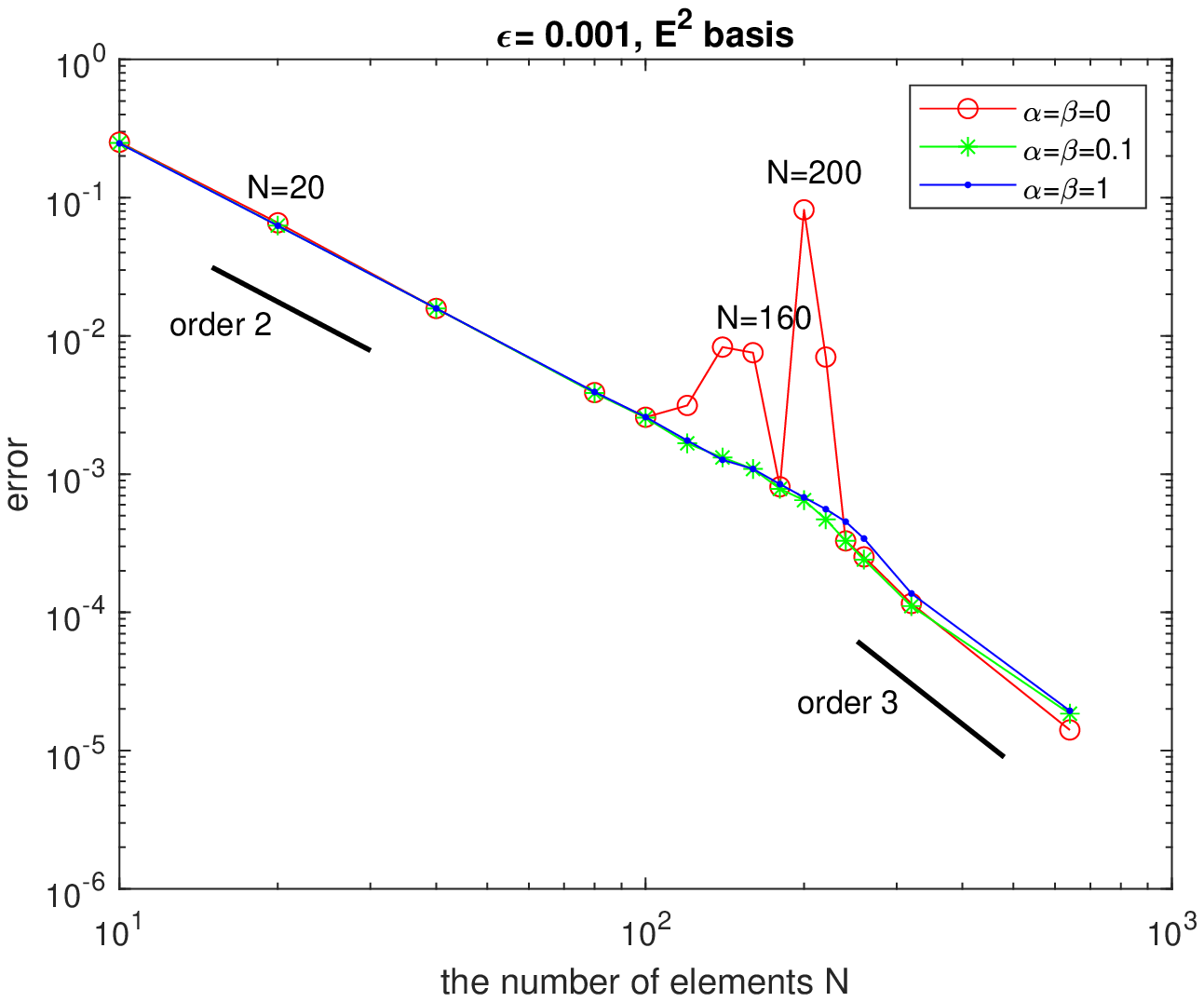}
    \includegraphics[width=2.6in,height=2.3in]{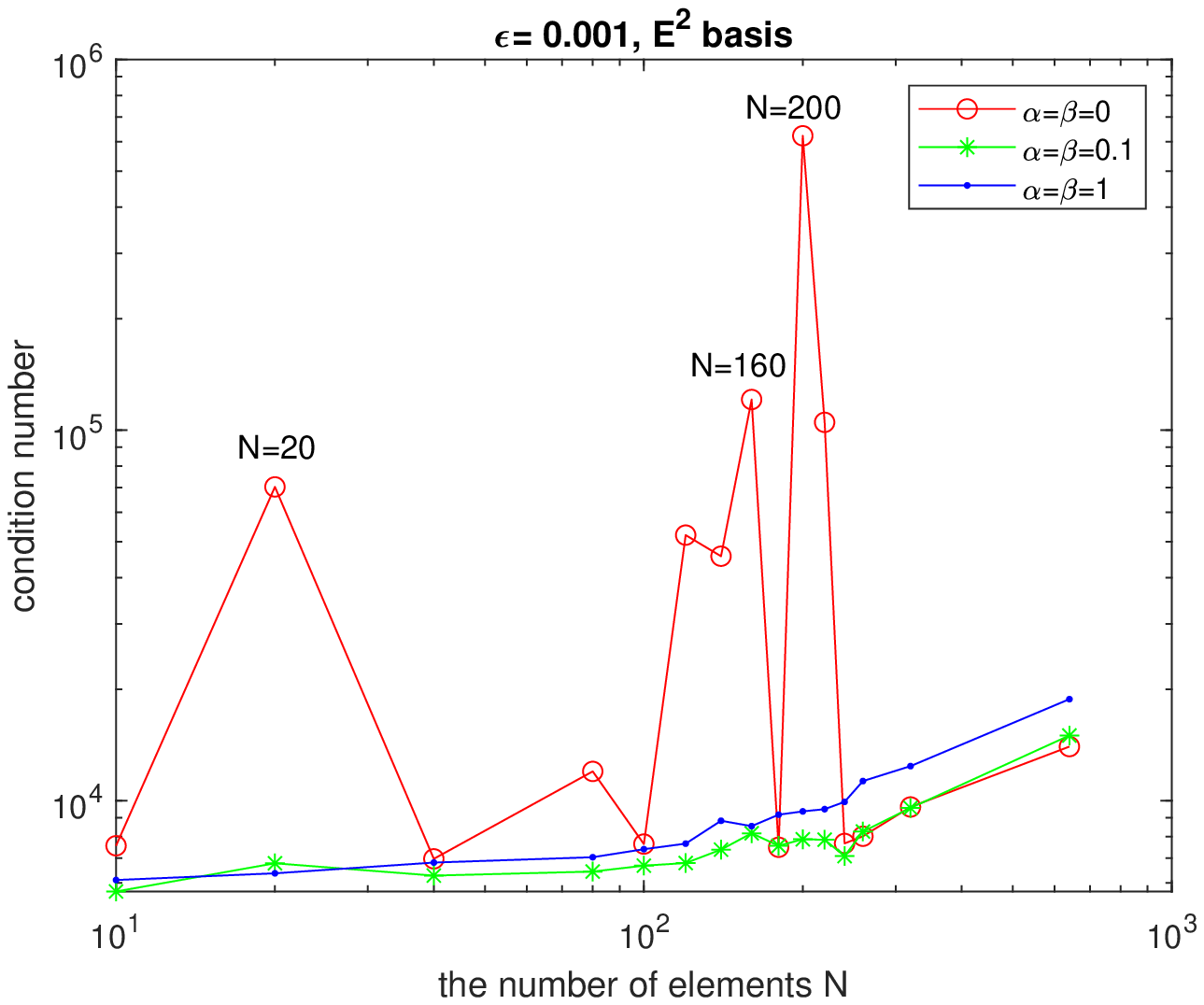}

\caption{Example \ref{ex:2}: Numerical results by multiscale DG ${E}^2$ for $\varepsilon=1\times 10^{-3}$. Left: $L^2$-errors of $u$. Right: condition numbers.}
\label{fig:1d3e2}
\end{figure}

\begin{figure}
  \centering
  \includegraphics[width=2.6in,height=2.3in]{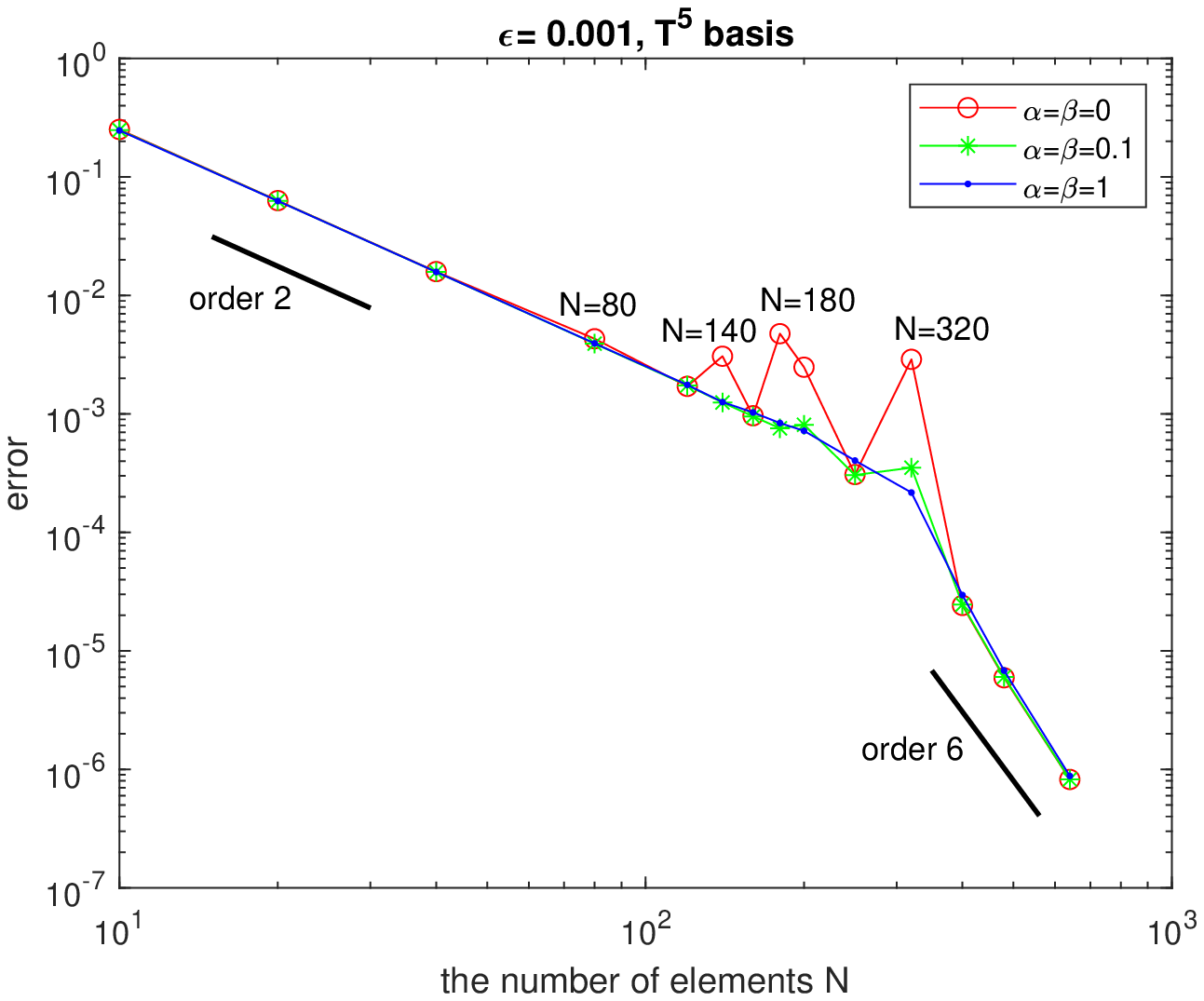}
    \includegraphics[width=2.6in,height=2.3in]{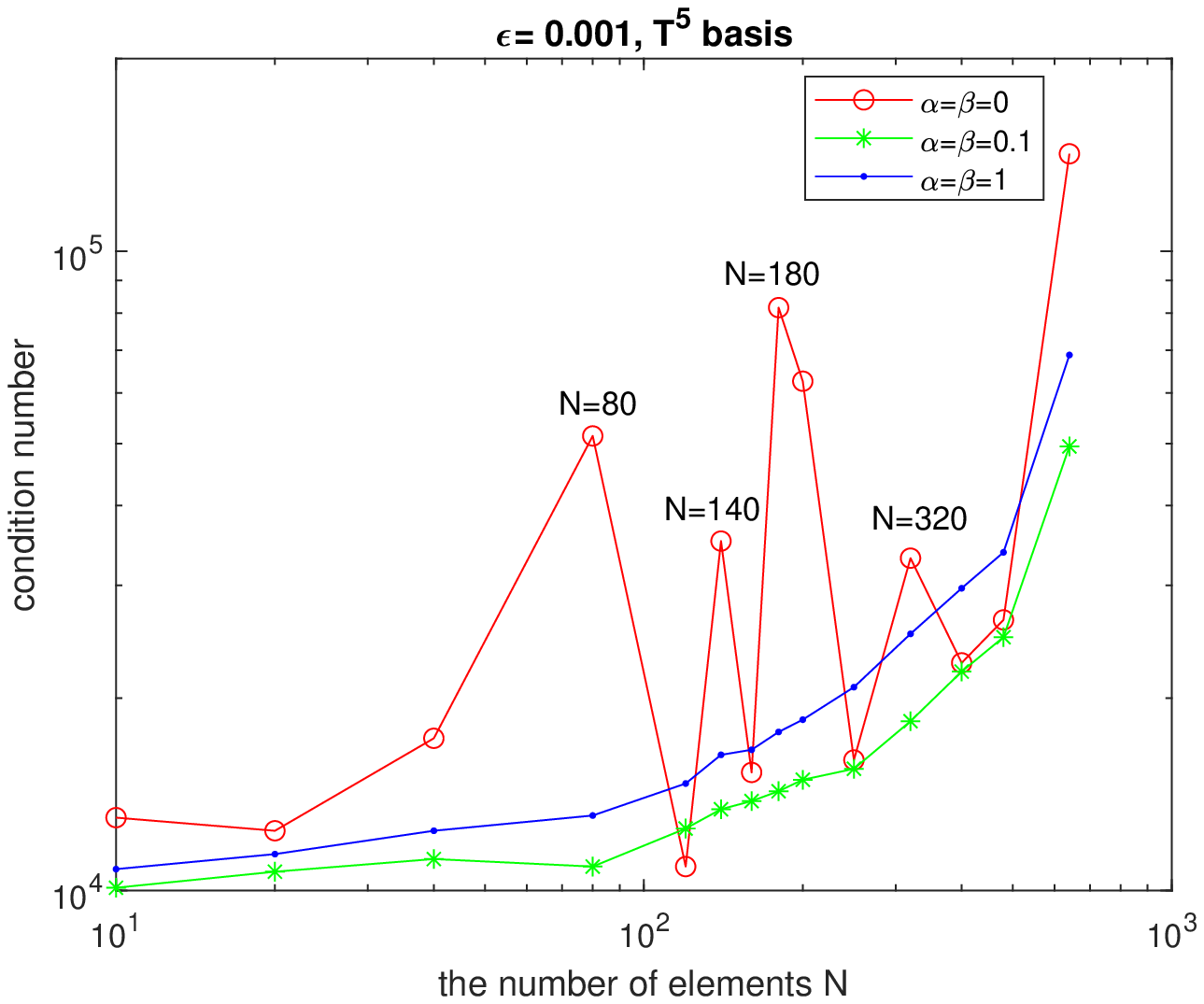}

\caption{Example \ref{ex:2}: Numerical results by multiscale DG ${T}^5$ for $\varepsilon=1\times 10^{-3}$. Left: $L^2$-errors of $u$. Right: condition numbers.}
\label{fig:1d3t5}

\end{figure}

When $\varepsilon$ reduces to $0.001$, % (i.e. $1/5$ of $0.005$, 
which is $1/5$ of $0.005$,
we expect the resonance location $N$ to be 5 times of that for $\veps=0.005$.
%the resonance location of number of element  is expected to  increase 5 times.)
%\bo{the number of elements $N$, where the resonance occurs, is expected to  increase by 5 times.} 
As we see in the error plot in Fig. \ref{fig:1d3e1}, the multiscale DG ${E}^1$ has resonance around $N=140$, which is about 5 times of previous $N=30$ in Fig. \ref{fig:5d3e1}. 
The condition numbers show spikes at  
the same locations.
 We can see similar phenomenon for multiscale DG with ${E}^2$, ${E}^3$,
 ${T}^3$ and ${T}^5$ as well; see Fig. \ref{fig:1d3e2}-\ref{fig:1d3t5}. 
%\bo{REMOVE (}Multiscale DG ${E}^2$ has resonance around $N=30$ and 45 for  $\varepsilon=0.005 $ and 160 and 200 for  $\varepsilon=0.001 $ (Fig.  \ref{fig:1d3e2}). 
%Multiscale DG ${E}^3$ has resonance around $N=25$ for  $\varepsilon=0.005 $ and 140 for  $\varepsilon=0.001 $ (Fig.  \ref{fig:1d3e3}).
%Multiscale DG ${T}^3$ has resonance from $N=40$ for  $\varepsilon=0.005 $ and from 140 to 320 for  $\varepsilon=0.001 $ (Fig.  \ref{fig:1d3t3}).
%Multiscale DG ${T}^5$ has resonance around $N=40$ and 70 for  $\varepsilon=0.005 $ and from 140 to 320 for  $\varepsilon=0.001 $ (Fig.  \ref{fig:1d3t5}).)

%Although  ${E}^2$,  ${E}^3$,  ${T}^3$ and ${T}^5$ can not get higher order \bo{approximations} than $E^1$ when $h\gtrsim\varepsilon$, they do have a smaller magnitude of errors than $E^1$.
%When $h\lesssim\varepsilon$, all of them can obtain the optimal order of accuracy, i.e.  $E^{p+2}$ has a $(p+3)$th order and  $T^{2p+1}$ has a $(2p+2)$th order.  In particular, we observe \bo{the} 3rd order for ${E}^2$, \bo{the} 4th order for ${E}^3$ and  ${T}^3$, and \bo{the} 6th order for ${T}^5$.
%

%-------------------------

\begin{figure}
  \centering
  \includegraphics[width=2.6in,height=2.3in]{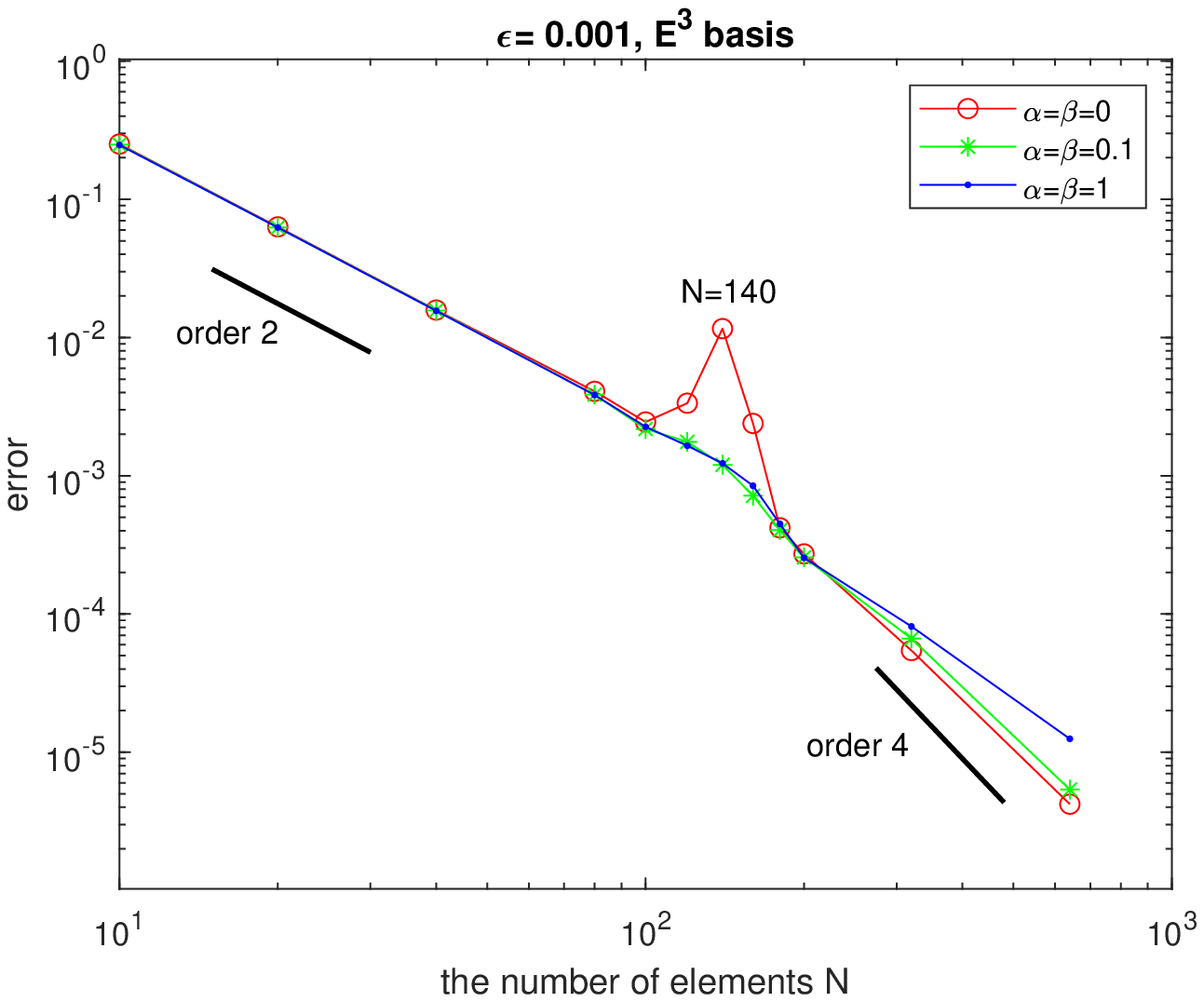}
    \includegraphics[width=2.6in,height=2.3in]{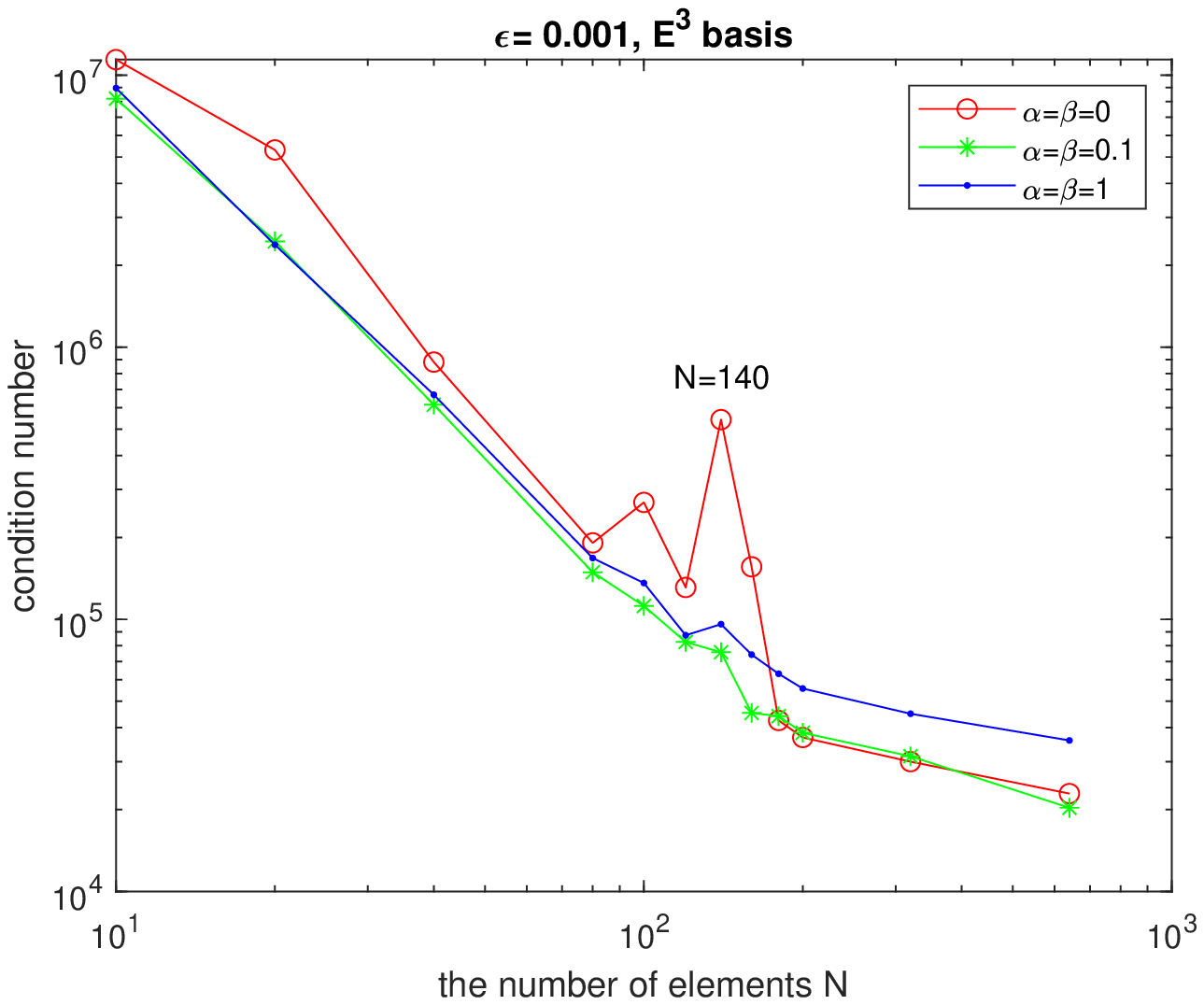}

\caption{Example \ref{ex:2}: Numerical results by multiscale DG ${E}^3$ for $\varepsilon=1\times 10^{-3}$. Left: $L^2$-errors of $u$. Right: condition numbers.}
%\label{fig:1d3e3}
\end{figure}

\begin{figure}
  \centering
  \includegraphics[width=2.6in,height=2.3in]{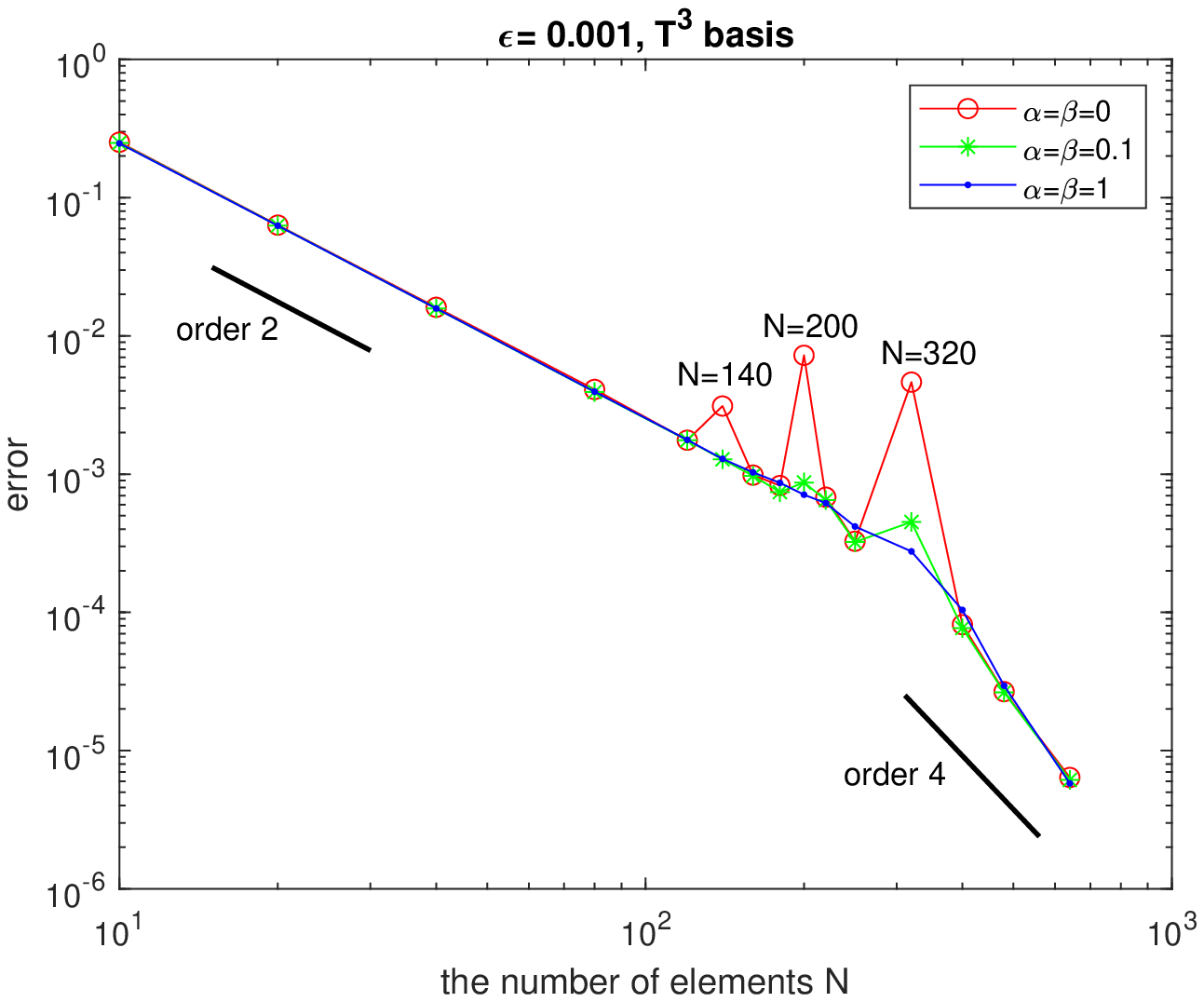}
    \includegraphics[width=2.6in,height=2.3in]{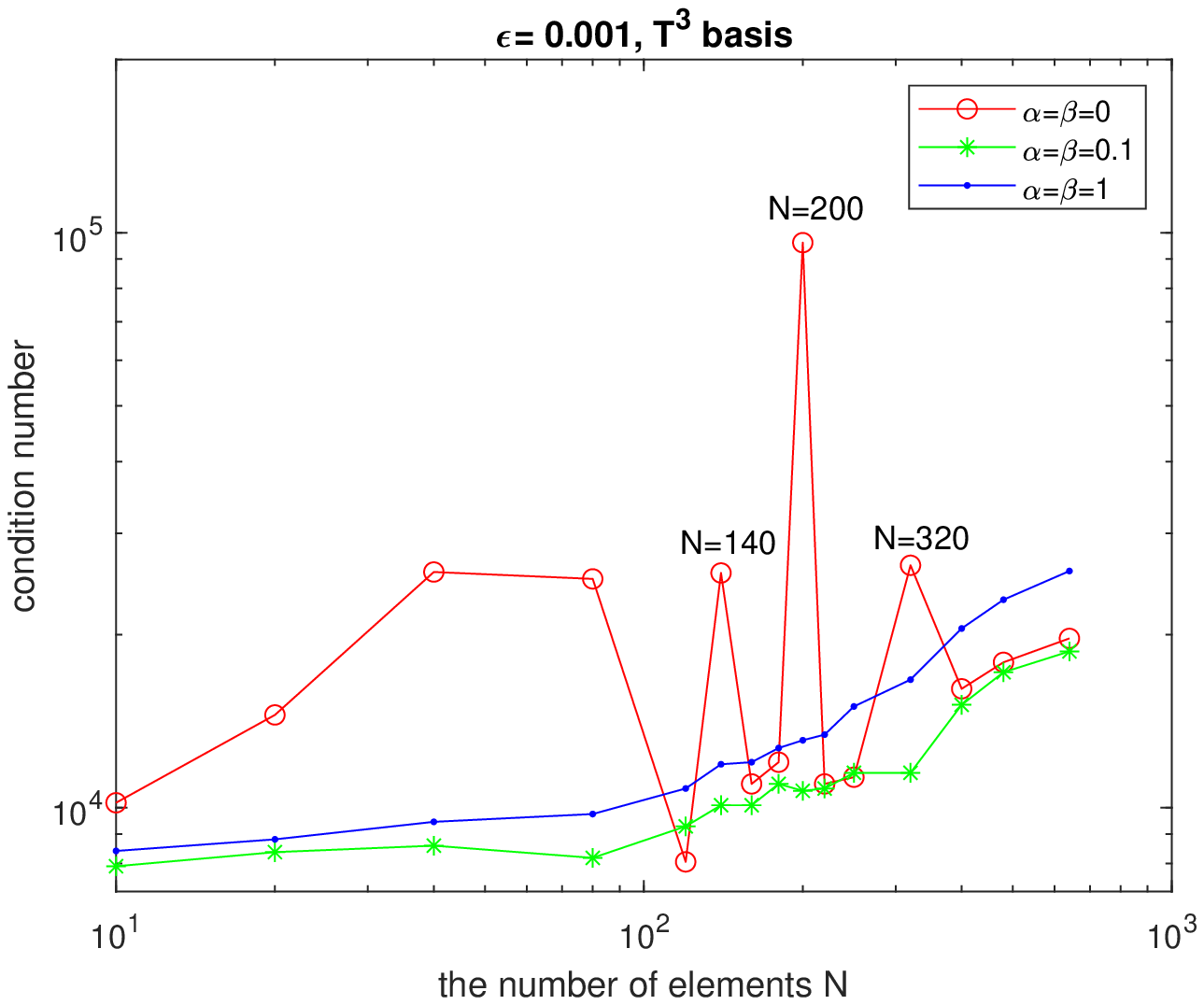}

\caption{Example \ref{ex:2}: Numerical results by multiscale DG ${T}^3$ for $\varepsilon=1\times 10^{-3}$. Left: $L^2$-errors of $u$. Right: condition numbers.}
%\label{fig:1d3t3}
\end{figure}

\section{Concluding remarks}

%(From our numerical experiments, we can conclude the following findings:)
In this paper, we carry out numerical experiments and have the following findings:

1. We observe resonance errors  in multiscale DG  with 
$\alpha=\beta=0$ for all  five
 different spaces   ${E}^1(=T^1)$,  ${E}^2$,  ${E}^3$,  ${T}^3$ and ${T}^5$. The multiscale DG with positive penalty parameters have almost no resonance errors. Thus %the
  taking positive penalty parameters significantly  helps 
 reduce the resonance errors.

2. In general, the resonance errors occur around $h\simeq  \varepsilon$. For different function spaces $E^p$ or $T^{2p-1}$,  the resonance errors may occur at different locations. Some may be around one mesh point, some may be in a small region. But the  graph of corresponding condition numbers has  spikes at the same locations.

3. Multiscale DG methods with positive penalty parameters show second order convergence when  $h \gtrsim  \varepsilon$ and
 optimal  order convergence when  $h \lesssim  \varepsilon$. Multiscale DG methods with $\alpha=\beta=0$ behave similiarly when the mesh size $h$ is away from the resonance location $h\simeq \veps$.

In the future work, we would like to investigate the resonance errors of  multiscale discontinuous Galerkin methods for two-dimensional Schr\"{o}dinger equation in \cite{BW22,BW22proceeding}. 
%Although no resonance errors were observed for zero-penalty multiscale DG for 2D Schr\"{o}dinger equation in \cite{BW22,BW22proceeding}, we believe the resonance errors exist. 
%We did not see it probably because the errors contain both $x$-direction  and $y$-direction. Our guess is that we won't see resonance errors unless the error from $x$-direction dominates. 

\section{Acknowledgments}
The authors consent to comply with all the Publication Ethical Standards. The research of the first author is supported
by National Science Foundation grant DMS-1818998.

\section{Conflict of interest}
The authors declare that they have no conflict of interest.

\end{document}